\documentclass[%
 reprint,
 amsmath,amssymb,
 aps,pre
]{revtex4-2}

\usepackage{tikz-cd} 
\usepackage{graphicx}
\usepackage{dcolumn}
\usepackage{bm}
\usepackage{amsthm,amssymb,mathtools}
\usetikzlibrary{arrows.meta,calc}
\usetikzlibrary{backgrounds}
\usepackage{subcaption}

\usepackage{booktabs}

\usepackage[shortlabels]{enumitem}

\graphicspath{{New_Figs/}}

\begin{document}

\preprint{APS/LK19830E}

\title{
Functorial invariants for chaos topology from data
}

\author{Denisse Sciamarella}
\email{denisse.sciamarella@cnrs.fr}
\affiliation{CNRS – Centre National de la Recherche Scientifique, 75016 Paris, France}
\affiliation{CONICET ‒ Universidad de Buenos Aires. Centro de Investigaciones del Mar y la Atmósfera (CIMA)}
\affiliation{Institut Franco-Argentin d'Études sur le Climat et ses Impacts (CNRS – IRD – CONICET – UBA) \\(IRL 3351 IFAECI), C1428EGA Ciudad Autónoma de Buenos Aires, Argentina}

\date{\today}

\begin{abstract}
The templex is a topological object bridging homologies and templates for chaotic dynamics. This article places the templex within category theory, introducing a directed path algebra, an edge operator on directed paths, and an equivalence relation for directed cycles that is distinct from directed homologies. The resulting functorial invariants are of two kinds: abelian-group invariants, namely the homology groups, and semigroup invariants, namely the generatex semigroups. These invariants are separable through forgetful functors and constitute a robust framework for identifying tipping points, disambiguating physical mechanisms, and benchmarking data-driven models against observations or simulations. The formulation sets forth a non-metric criterion for chaos from finite-time data and reveals that the concatenable nature of Topological Modes of Variability is a direct consequence of the semigroup structure of the directed path algebra. Two applications are presented: an experimental speech signal and a climatic numerical simulation. 
\end{abstract}

\maketitle

\section{Introduction}

Complex systems are governed by nonlinear dynamics exhibiting regimes, phases, transitions, tipping points, or extreme events. Despite the increasing abundance of data from observations and simulations, a persistent gap remains between what can be measured and the understanding of the main mechanisms governing the observed dynamics. Current approaches rely on decompositions based on linear or statistical methods which struggle to capture the structural constraints imposed by the dynamics itself. This limitation is widely acknowledged in recent reviews, which emphasize the need for frameworks capable of bridging data-driven approaches with fundamental principles \cite{brunton2020machine}.

Closing this gap requires topological invariants, without which the structure of a chaotic flow cannot be robustly reconstructed \cite{gilmore1998topological, ghil2023dynamical}.  The Lorenz and Rössler chaotic attractors, for instance, exhibit nearly identical fractal dimensions, rendering them metrically indistinguishable under observational noise \cite{Let10b}, in spite of their distinct features.  Chaos topology addresses this issue: it provides a framework to link data to the principles governing trajectory organization in phase space. Under this paradigm, the topological characterization of flows becomes a central problem \cite{gilmore2012topology}, relevant when validating reduced-order models, comparing observational sources, or assessing the consistency of models extracted via data-driven identification methods as well as physics-informed machine learning \cite{karniadakis2021physics}.

Understanding chaotic dynamics from trajectory data has traditionally required two steps: reconstructing orbits in phase space, and then describing how they wind around one another using braids, knots, or linking numbers \cite{birman1983knotted, gilmore1998topological}. This information can be condensed in a so-called template or knot-holder. In a template, knotted orbits cycle around a joining chart where subsets of trajectories with different paths in phase space, called \emph{strips}, merge together and share a common path before splitting apart again. This approach is effective in three dimensions for low-noise datasets, but becomes impractical as soon as orbits cannot be reliably recovered from point clouds, or in higher-dimensional dynamics, where knots generically unknot. Several attempts have explored how to extend templates beyond three dimensions, and move away from knot-based descriptions \cite{lefranc2006alternative,mangiarotti2014modelisation}.

A different approach bypasses orbit reconstruction altogether by working directly with the structure that supports a point set in phase space. This approach led to the construction of a cell complex in the sense of algebraic topology \cite{sciamarella1999topological,sciamarella2001unveiling}. A cell complex is a scaffold composed of cells of layered dimensions, built in such a way that the topological invariants, i.e. properties that are preserved under homeomorphisms, can be extracted independently of the particular distribution or number of cells involved \cite{kinsey2012topology}. Constructing cell complexes from point sets can be done following different rules, which will not necessarily yield a topologically faithful representation of the manifold supporting the underlying dynamics. Branched Manifold Analysis through Homologies assumes that the point set lies on a manifold which may have boundaries and branches, and builds a particular type of cell complex, now called BraMAH complex \cite{charo2020topology,ghil2023dynamical}, especially designed for this task. This framework allows one to compute homology groups (holes of different dimensions) of the supporting manifold, retaining generators and orientability properties. One can thus distinguish between spaces such as Möbius bands and cylinders, tori, Klein bottles, and pseudo-manifolds, such as branched manifolds.

Classical algebraic topology, however, describes only the structure of a space, without accounting for how that structure is actually visited by a flow. The domain that deals with the topology of directed spaces was born in the 1990s and is called directed algebraic topology \cite{grandis2003directed,goubault2003some, gaucher2003homotopical, grandis2009directed}. To incorporate the directed organization of a flow in phase space, one must then move from traditional algebraic topology to directed algebraic topology. In parallel, when the goal is to describe chaotic attractors, one should move from a description in terms of BraMAH complexes to a description in terms of templexes \cite{charo2022templex, sciamarella2024new, mosto2024templex}, which encompasses the latter. A templex is a BraMAH cell complex endowed with a directed graph (digraph) whose nodes are the locally highest-dimensional cells (top-cells) and whose edges are the flow-compatible (directed) connections between these nodes. The combination of both in a single object makes this digraph fundamentally different from those used in network analysis or clustering tools \cite{zou2019complex}. Because the digraph is tied to a cell complex, this structure allows for the extraction of invariants regardless of the specific segmentation of the point cloud, thus canceling out the scaffold's extraneous information. The templex properties thus include the directionless invariants of classical topology mentioned above as well as the directed ones, called \emph{stripexes}, since they extend the concept of strips in a template to higher-dimensional dynamics; for instance, the four strips of the classical Lorenz template correspond to four stripexes in the Lorenz templex. As we shall see, these directed invariants involve an equivalence relation that does not rely on topological voids. This distinction is central: while the undirected invariants describe the presence of holes, torsions or junction loci, the directed invariants in a templex provide information about the specific type of chaotic organization at play.

Historically, the definition of cell complexes \cite{whitehead1949combinatorial} and their associated homology preceded the advent of modern computational topology by several decades. This article follows the inverse path, providing a formal foundation for computational definitions of a templex and its properties \cite{sciamarella2023code}.  This is achieved by leveraging the language of category theory \cite{EilenbergMacLane1945,MacLane1971}, which emerged in the 1940s as a technical tool in algebraic topology. A category is defined by a collection of objects (e.g., spaces, groups) and arrows that represent morphisms between them (e.g., continuous maps, homomorphisms). A central notion in category theory is that of a \emph{functor}, formalizing the idea that a construction in mathematics should respect the relations between objects, not just the objects themselves. For instance, the homology functor translates topological spaces into abelian groups, and continuous maps into group homomorphisms, preserving composition.

This article introduces several elements: a new type of functor, an operator acting on finite formal concatenations of directed edges, and the quotient of directed cycles modulo the image of this operator. The algebraic structure that emerges for the directed invariants of a templex is that of a semigroup. Two applications are presented: a templex from an experimental speech signal is constructed for the first time, showing how the directed invariants enhance the homological description \cite{sciamarella1999topological}, and a recent application to non-autonomous climate dynamics is revisited, including the use of directed invariants to define the Topological Modes of Variability (TMVs). Previously introduced and extracted algorithmically \cite{charo2025topological}, TMVs are here given a categorical formulation and a structural explanation for their concatenable character.  The article concludes with an outline of how the functorial invariants can be interpreted in terms of the mechanisms organizing chaos in phase space, together with topological perspectives on the two applications and, more generally, on nonlinear time series analysis. For completeness, the appendices provide a summary of standard notions from homology theory, and illustrate how the functorial invariants emerge from the Rössler and Lorenz attractors \cite{lorenz1963lorenz,Ros76c}. 

\section{The templex: definitions}
\label{sec:templex}

This section introduces mathematical definitions of the templex that generalize earlier constructions and place previously algorithmic procedures within a rigorous framework, starting with an intrinsic definition of a BraMAH complex and then incorporating flow-directionality.

Cell complexes can be of different types depending on the construction rules used to build them from a point cloud. A BraMAH complex is designed so that the topology of a (branched) manifold underlying a point cloud in phase space can be faithfully described by a cell complex whose dimension coincides with that of the underlying manifold \cite{sciamarella1999topological, sciamarella2001unveiling,charo2020topology}. Here, a branched manifold is understood as a space that is locally a manifold almost everywhere~\cite{gilmore1998topological}. In Ref.~\cite{sciamarella2024new}, a BraMAH complex of dimension $h$ is defined as a cell complex constructed from a point cloud embedded in an $n$-dimensional phase space, with cells of dimension $k\le h$, such that: (i) the $0$-cells are a sparse subset of the original cloud of points in phase space; (ii) $h=d$, where $d\le n$ is the intrinsic dimension of the underlying branched manifold being sampled; and (iii) each $h$-cell is such that $h$ of the singular values that describe the distribution of the points around the barycenter of the cell scale linearly with the number of points in the cell~\cite{sciamarella1999topological, sciamarella2001unveiling}.  

A BraMAH complex can be defined independently of algorithmic implementations as a finite cell complex~$K$ whose geometric realization $|K|$ is a manifold or a branched manifold. A cell is said to be top-dimensional if it is not a face of any cell of strictly higher dimension in the complex, even if its dimension is not maximal globally. A \emph{junction locus} is a subcomplex where three or more top-dimensional cells meet along cells of lower dimension, assuming a minimal local cell structure; see Ref.~\cite{charo2022templex} for details. These are the loci where a branched manifold fails to be a manifold.

Computing homology groups for a BraMAH complex provides an algebraic translation of the topological properties of the structure of a possibly branched manifold sampled by a point cloud in phase space. The zeroth homology group detects the number of connected components of the sampled support, indicating whether the system explores a single connected region of phase space or several disconnected ones. Homology groups of degree $k \ge 1$ capture higher-dimensional features, corresponding to regions of phase space avoided by the data and translating into voids of the support~\cite{sciamarella1999topological, sciamarella2001unveiling, charo2020topology, charo2021topological}.  Unlike constructions commonly used in persistent homology~\cite{edelsbrunnerHarer2010}, the structure of a BraMAH complex does not depend on any filtration parameter. Top-cells can gather large subsets of points sampling the local topological structure of the manifold. This facilitates the computation of homology groups with integer coefficients, alongside the identification of extra structural features such as junction loci.

While classical algebraic topology captures the global structure of the set sampled in phase space, it is insensitive to temporal ordering. Let $K$ be a BraMAH complex, and let $\mathcal{G} = (N,E)$ be a directed graph, or {\em digraph}, whose nodes $N$ correspond one-to-one with the top-cells of $K$, encoding the flow-compatible connections between them. The pair $\mathcal{T} = (K,\mathcal{G})$ defines a \emph{templex}. A \emph{subtemplex} of $\mathcal{T}$ is a pair $(K',\mathcal{G}')$ where $K'\subset K$ is a subcomplex and $\mathcal{G}'$ is the subgraph of $\mathcal{G}$ induced on the top-cells of $K'$.  Once the digraph $\mathcal{G}$ is taken into account, the junction locus, previously defined as a subcomplex of $K$, naturally induces a subtemplex of $\mathcal{T}$ whose top-cells can be identified as \emph{ingoing} or \emph{outgoing} according to the direction of the flow across the lower-dimensional cell forming the junction. An \emph{ingoing cell} (resp.\ \emph{outgoing cell}) may therefore be referred to equivalently as an \emph{ingoing node} (resp.\ \emph{outgoing node}). Following the terminology introduced in template theory and used throughout previous works \cite{charo2022templex, charo2025topological}, junction loci with a single ingoing node are called \emph{splitting loci}, while those with a single outgoing node are called \emph{joining loci}.

A \emph{generatex} is a subtemplex whose subgraph is a directed cycle of $\mathcal{G}$. It is said to be of order $p$ with $p \in \mathbb{N}$, $p \ge 1$, if the cycle has $p$ distinct ingoing nodes. A generatex of order $1$ is a stripex and a generatex of order $p>1$ can be decomposed into $p$ stripexes. Stripexes provide the direct bridge back to strips in the template tradition. By construction, any templex endowed with at least one joining locus necessarily contains intersections between its generatexes. For any subset of indices $\{i_1,\dots,i_m\}$ with $m\ge 2$, we define a \emph{bond} as $B_{i_1\cdots i_m} = G_{i_1}\cap\cdots\cap G_{i_m}$, whenever this intersection is non-empty. The index set records precisely which generatexes share that portion of the dynamics, and the \emph{valence} of the bond is $m$. The term \emph{valence} refers to the analogy with valence electrons in chemistry, to emphasize that bonds encode how the fundamental dynamical units combine through the gluing of generatexes, in the same spirit as introduced in Ref.~\cite{mosto2024templex}. The union of generatexes can be represented as a directed multigraph, i.e. a graph with multiple directed edges between the same two nodes.  In this multigraph, the valence of a given bond coincides exactly with the number of parallel edges along that subtemplex. Bonds with different index sets are distinct, even if they happen to share edges, as they encode different gluing relationships among the generatexes.

\section{Algebra of directed paths}

Following classical homology, this section develops an algebra of directed paths for a templex, conceived as the analogue of the algebra of chains for a cell complex. Let $\mathcal{T} = (K, \mathcal{G})$ be a templex and let $C_{\mathcal{G}}(\mathcal{T})$ denote the set of all finite formal concatenations of directed edges of $\mathcal{G}$. Concatenation defines an associative product, so that for any three paths $p,q,r$,
\[
(p\,q)\,r = p\,(q\,r).
\]
As direction matters in a templex, $C_{\mathcal{G}}(\mathcal{T})$ has the structure of a semigroup. Among the directed edges of a templex there are special ones: those connecting several ingoing nodes to a single outgoing node implicitly contain a lower-dimensional cell serving as frontier between those cells or nodes that are exclusive to a single generatex and those that belong to several generatexes. This lower-dimensional cell constitutes the best location for Poincaré sections to compute strips with first return maps in template theory. These key directed edges are hereafter called \emph{Poincaré edges}. A Poincaré edge is denoted symbolically by a pair of ingoing and outgoing node labels,
\[
\langle \mathrm{ingoing~node~label}\mid 
\mathrm{outgoing~node~label}\rangle,
\]
and we denote by $\bar{E}$ the set of all Poincaré edges of the templex $\mathcal{T}$.  Notice that each joining locus gives rise to at least two Poincaré edges, reflecting the fact that the same outgoing top-cell/node is reached from two or more adjacent top-cells/nodes.  The terminology is chosen to emphasize the role of these directed edges as transverse dynamical interfaces where distinct directed paths merge to become indistinguishable. The \emph{Poincaré-edge operator} is defined as
\[
\mathcal{P} : C_{\mathcal{G}}(\mathcal{T}) \longrightarrow \bar{E}^*,
\]
where $\bar{E}^*$ denotes the set of all finite ordered sequences of Poincaré edges in $\mathcal{T}$.  Given a directed path $p$, the operator $\mathcal{P}(p)$ returns the ordered list of joining directed edges along $p$, in the order induced by the path, discarding all other directed edges. Thus, operationally, $\mathcal{P}$ contracts a directed path onto the Poincaré edges it traverses. For directed cycles, this sequence is understood up to cyclic rotation. The \emph{path templex} associated with a templex $\mathcal{T}$ is the algebraic structure
\[
T_\bullet(\mathcal{T}) = \bigl( (C_\bullet(K),\,\partial),\,
(C_{\mathcal{G}}(\mathcal{T}),\,\mathcal{P}) \bigr),
\]
where $C_\bullet(K)$ denotes the chain complex associated with a BraMAH complex $K$, and $\partial$ the boundary operator (Appendix~\ref{app:homology}). The notion of path templex introduced here should not be confused with that of a directed chain complex, nor with the dihomology framework developed in directed algebraic topology \cite{grandis2003directed,goubault2003some, gaucher2003homotopical, grandis2009directed}. In the present construction, the underlying chain complex $C_\bullet(K)$ remains entirely classical: no homological theory is modified or replaced. Causality is encoded exclusively by directed paths and by the Poincaré-edge operator acting on them. The algebraization preserves a separation between the structure visited by the flow and the flow upon it, without decoupling them.

Let $\mathsf{Z}_{\mathcal G}$ denote the set of directed cycles, i.e., of closed directed paths, in the digraph $\mathcal G$. Two directed cycles in $\mathcal G$ are said to be equivalent if and only if they have the same image under the Poincaré-edge operator. This equivalence relation induces the quotient
\[
\mathrm{Gen}(\mathcal{T}) = \mathsf{Z}_{\mathcal G} / {\sim},
\]
whose elements are the \emph{generatex classes}. By a standard abuse of notation, we will use the term \emph{generatex} and the symbol $G_i$ to denote both a specific directed cycle (the subtemplex) and its corresponding causal equivalence class in the semigroup.  This construction mirrors the classical algebraic-topological pipeline: directed paths play the role of chains, the Poincaré-edge operator extracts boundary-like causal information, and generatex classes arise as a quotient of directed cycles by this operator. Notice that regular or quasiperiodic dynamics have no junction locus, hence no Poincaré edges: all directed cycles collapse into a single trivial equivalence class under $\mathcal{P}$, and the generatex semigroup reduces to the trivial semigroup.  Conversely, a junction locus structurally forces the merging of at least two non-equivalent directed paths, yielding a non-trivial generatex semigroup.

Let $\mathcal{T}_G = (K_G, \mathcal{G}_G)$ be the subtemplex defined by the union of all directed cycles belonging to a specific generatex class $G$. Assigning a consistent local orientation to each top-cell of $K_G$ by propagating it across shared boundaries, the orientability chain is the sum of all their signed boundaries. A class $G$ is \emph{orientation-reversing} if any cell of codimension 1 appears in this chain with a coefficient of absolute value greater than~$1$; otherwise it is \emph{orientation-preserving}.

\section{Functorial invariants}
\label{sec:functorial}

The constructions introduced in the previous sections reveal two complementary levels in the description of dynamics: a structural one, captured by the topology of the underlying space, and a directional one, reflecting how this structure is effectively visited by the flow.  This section addresses how both directionless and directed invariants can be related through a unified categorical perspective.  

Let $\mathbf{Top}$ be the category of topological spaces and continuous maps, and $\mathbf{Ab}$ that of abelian groups and group homomorphisms, while $\mathbf{Ch}$ has graded abelian groups endowed with a boundary operator as objects, and degree-preserving morphisms commuting with the boundary operator as arrows. The classical functorial chain $\mathbf{Top} \longrightarrow \mathbf{Ch} \longrightarrow \mathbf{Ab}$ is represented by the following commutative diagram, where the homology functor maps each chain complex to its homology groups:
\[
\begin{tikzcd}[column sep=huge, row sep=large]
X \arrow[r] \arrow[d] & C_\bullet(X) \arrow[r] \arrow[d] & H_k(X) \arrow[d]\\
Y \arrow[r] & C_\bullet(Y) \arrow[r] & H_k(Y)
\end{tikzcd}
\]
In directed algebraic topology, the category $\mathbf{dTop}$ consists of directed spaces and direction-preserving maps. In this setting, classical cell complexes and chain algebras are replaced by templexes and algebras of directed paths. Let us denote by $\mathbf{Pa}$ the category whose objects are directed path semigroups endowed with the Poincaré-edge operator, and whose morphisms are semigroup homomorphisms that send Poincaré edges to Poincaré edges.

A forgetful functor maps an object to a structure obtained by discarding part of the object's data, while keeping the remaining structure and the induced morphisms. Two forgetful functors arise naturally in this setting. The first one, $U_{\mathcal{G}} : \mathbf{dTop} \longrightarrow \mathbf{Ch}$, forgets the digraph component and retains only the underlying chain complex, thereby recovering the classical functorial chain.  The second one, $U_{K} : \mathbf{dTop} \longrightarrow \mathbf{Pa}$, forgets the chain complex structure but preserves the directed path algebra together with the Poincaré-edge operator, which inherits properties derived from the prior construction of the BraMAH complex. These properties are thus translated into the non-graded category $\mathbf{Pa}$: causal direction does not propagate to lower-dimensional cells, but remains confined to transitions between top-cells.

The categories $\mathbf{Ch}$ and $\mathbf{Pa}$ are scaffold-dependent. Functorial invariants must condense the scaffold-independent properties. The generatex functor, $F_{\mathrm{Gen}} : \mathbf{Pa} \longrightarrow \mathbf{Sem}$, maps each path structure to the semigroup generated by its generatex classes, $\mathbf{Sem}$ being the category of semigroups and semigroup homomorphisms. This functor plays for directed paths in a templex a role conceptually parallel to that played by the homology functor for chain complexes. The resulting composition, $F_{\mathrm{Gen}} \circ U_{K} : \mathbf{dTop} \longrightarrow \mathbf{Sem}$, defines a functorial construction parallel to the classical homology theory, but defining causal equivalence classes, and leading to a semigroup encoding its intrinsic directed invariants.  
\[
\begin{tikzcd}[column sep=huge, row sep=large]
X^\uparrow \arrow[r] \arrow[d, ""'] 
  & T_\bullet(X^\uparrow) \arrow[r] \arrow[d, ""'] 
  & (H_k(X),\mathrm{Gen}(X^\uparrow)) \arrow[d, ""'] \\
Y^\uparrow \arrow[r] 
  & T_\bullet(Y^\uparrow) \arrow[r] 
  & (H_k(Y),\mathrm{Gen}(Y^\uparrow))
\end{tikzcd}
\]
The new functorial chain provides a framework in which the classical directionless invariants, i.e., the homology groups, as well as the directed invariants, i.e., the generatex semigroups, coexist. While homology group generators are classes of chains encircling voids in a topological space, generatex semigroups are classes of directed cycles circling around joining loci.  The analogies between the two kinds of functorial invariants are summarized in Table~\ref{tab:ch_pa_comparison}. 

\begin{table}[!htbp]
\centering
\renewcommand{\arraystretch}{1.5}
\begin{tabular}{l | l l}
\hline
\textbf{Level} & \textbf{Directionless} & \textbf{Directed} \\
\hline
\begin{minipage}[t]{0.25\columnwidth}\raggedright Spaces\end{minipage} & 
\begin{minipage}[t]{0.30\columnwidth}\raggedright Topological spaces $X \in \mathbf{Top}$\end{minipage} & 
\begin{minipage}[t]{0.30\columnwidth}\raggedright Directed spaces $X^\uparrow \in \mathbf{dTop}$\end{minipage} \\

\begin{minipage}[t]{0.25\columnwidth}\raggedright Algebraic encoding\end{minipage} & 
\begin{minipage}[t]{0.30\columnwidth}\raggedright Chain complex algebra $(C_\bullet,\partial) \in \mathbf{Ch}$\end{minipage} & 
\begin{minipage}[t]{0.30\columnwidth}\raggedright Directed path algebra $(C_{\mathcal G},\mathcal P) \in \mathbf{Pa}$\end{minipage} \\

\begin{minipage}[t]{0.25\columnwidth}\raggedright Algebraic invariants\end{minipage} & 
\begin{minipage}[t]{0.30\columnwidth}\raggedright Homology groups $H_k(X) \in \mathbf{Ab}$\end{minipage} & 
\begin{minipage}[t]{0.30\columnwidth}\raggedright Generatex semigroups $\mathrm{Gen}(X^\uparrow) \in \mathbf{Sem}$\end{minipage} \\
\hline
\end{tabular}
\caption{Categorical correspondences.
}
\label{tab:ch_pa_comparison}
\end{table}

From a categorical perspective, a \emph{pushout} provides a canonical way of gluing two algebraic structures along a common substructure. Accordingly, bonds between generatex classes admit a natural pushout interpretation. Let $\mathrm{Gen}(X^\uparrow)\in\mathbf{Sem}$ denote the generatex semigroup associated with a directed space $X^\uparrow$. Under $F_{\mathrm{Gen}}$, the bond introduced in Section~\ref{sec:templex} as an intersection of directed cycles becomes an intersection of subsemigroups: given $G_1,G_2\le \mathrm{Gen}(X^\uparrow)$, their \emph{bond} is the common subsemigroup $B_{12}=G_1\cap G_2$.  The canonical inclusions $B_{12}\hookrightarrow G_1$ and $B_{12}\hookrightarrow G_2$ give rise to the pushout diagram
\[
\begin{tikzcd}[column sep=huge, row sep=huge]
B_{12} \arrow[r, hook, "i_1"] \arrow[d, hook', "i_2"]
  & G_1 \arrow[d, "j_1"] \\
G_2 \arrow[r, "j_2"']
  & G_1 \sqcup_{B_{12}} G_2
\end{tikzcd}
\qquad (\text{in }\mathbf{Sem}),
\]
which formalizes the gluing of generatex subsemigroups along their bond. 

Within the functorial framework developed here, stripexes do not require a separate categorical treatment. They are elements of the generatex semigroup $\mathrm{Gen}(\mathcal{T})$, and the generatex functor $F_{\mathrm{Gen}}$ maps them into $\mathbf{Sem}$.

\section{Decomposition into Topological Modes of Variability}
\label{sec:TMV}

The functorial invariants introduced in Section~\ref{sec:functorial}  are, by construction, metric-free: they describe equivalence classes of directed cycles, not the specific trajectories that realize them. A  complementary, metric-level description is obtained by identifying the segments that correspond to each  generatex class. This identification defines the \emph{Topological Modes of Variability} (TMVs) introduced in Ref.~\cite{charo2025topological}: a TMV consists of all trajectory segments that traverse the cells of a given generatex.

Let $\mathcal{T} = (K,\mathcal{G})$ be a templex and let $\mathrm{Gen}(\mathcal{T}) = \{G_1, \dots, G_n\}$ be its set of generatex classes. A trajectory $\phi : [0,T] \to |K|$ induces a sequence of top-cell visits, which can be read off from the digraph $\mathcal{G}$. We define a \emph{generatex visitation label} of $\phi$ as a map
\[
\chi : [0,T] \longrightarrow \{G_1,\dots,G_n\},
\]
assigning to each instant $t$ the generatex label carried by the trajectory at that time.  For generatex classes of order $p>1$, the trajectory may cross a joining locus and receive a new label before completing a full cycle, since each such crossing presents a new switching opportunity.  The labeling $\chi(t)$ is therefore updated at each crossing of a joining locus (the lower-dimensional cell in $K$ where the branched manifold fails to be a manifold) and remains constant between successive crossings. When the trajectory enters a bond, it traverses cells that belong to several generatex classes simultaneously. The labeling $\chi$ does not become undefined or ambiguous in this region: the bond interval inherits the label of the specific generatex class the trajectory eventually resolves into.

The generatex visitation label can be used to partition the time interval $[0,T]$ of (a set of) time series into a sequence of subintervals. The restriction of the trajectory to a subinterval labeled $G_i$ is one \emph{instance} of TMV-$i$. The full set of such subintervals, ordered in time, constitutes the \emph{TMV decomposition} of the trajectory:
\[
[0,T] = I_1 \cup I_2 \cup \cdots \cup I_N, 
\quad \chi\vert_{I_k} = G_{i_k},
\]
where $i_k \in \{1,\dots,n\}$ is the label of the generatex visited during $I_k$, and consecutive intervals are labeled by generatexes that are adjacent in the digraph $\mathcal{G}$. The semigroup structure of the directed path algebra, whose operation is concatenation, makes this decomposition \emph{concatenable}: the subintervals are placed end to end, not added. This stands in direct contrast with linear-based modal decompositions (Fourier, Empirical Orthogonal Functions, Dynamic Mode Decomposition, Singular Spectrum Analysis), where modes are superposed at every instant.

All the notions introduced so far are worked out explicitly for the Rössler and Lorenz attractors in Appendix~\ref{app:RL}.

\section{Applications}

Two case studies illustrate how the functorial invariants emerge from experimental \cite{sciamarella1999topological} and numerical \cite{charo2025topological} data. The first example is a speech signal previously used to show that the topological structure of flows reconstructed from short and noisy time series can be partially retrieved through homologies \cite{sciamarella1999topological}. The templex and its directed invariants are computed here for the first time. While the correspondence between the embedded trajectory and a three 1-hole structure cannot be taken as an indication of chaos, the directed invariants obtained through the generatex functor can, as shown below. The second example is a numerical simulation of a conceptual four-dimensional model~\citep{Pier.Ghil.2021} that captures key processes of the large-scale ocean circulation.  This system is revisited to show how functorial invariants can be handled in nonautonomous settings, where the TMVs can be used to detect the {\em topological tipping points} (TTPs) found in Ref.~\cite{charo2025topological}.

\subsection{A speech signal}

Vowel sounds result from a self-sustained oscillatory regime in which vocal-fold motion couples with the airflow driven by subglottal pressure. The vowel quality is further shaped by the supraglottal vocal tract, which acts as an acoustic resonator. The governing equations are inherently nonlinear, involve delayed feedback mechanisms, and are therefore known to produce complex and sometimes chaotic dynamical signatures.  

\begin{figure}[t]
\centering
(a)\includegraphics[width=0.7\columnwidth]{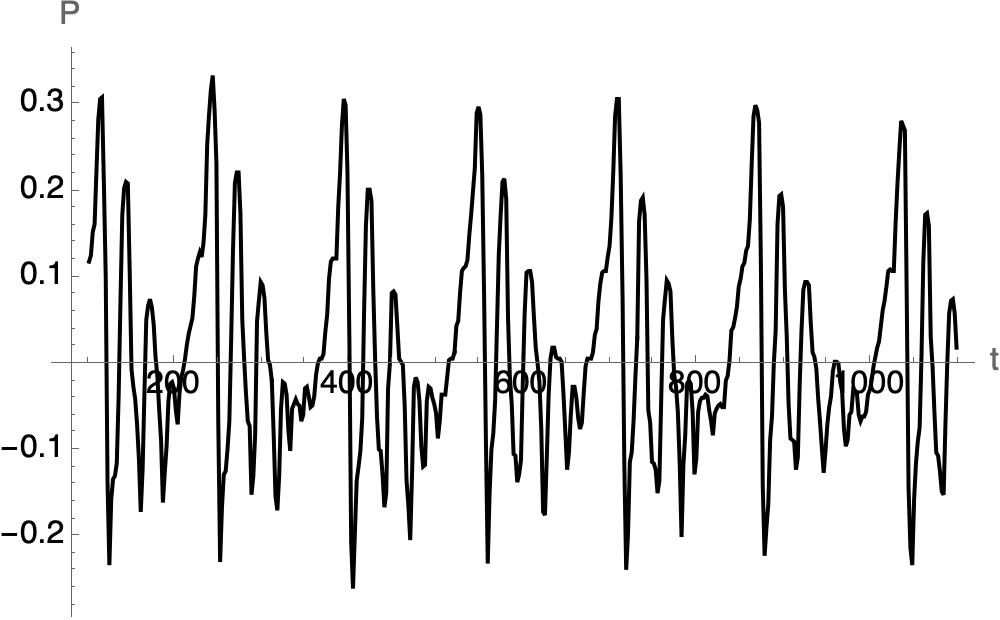}
(b)\includegraphics[width=0.7\columnwidth]{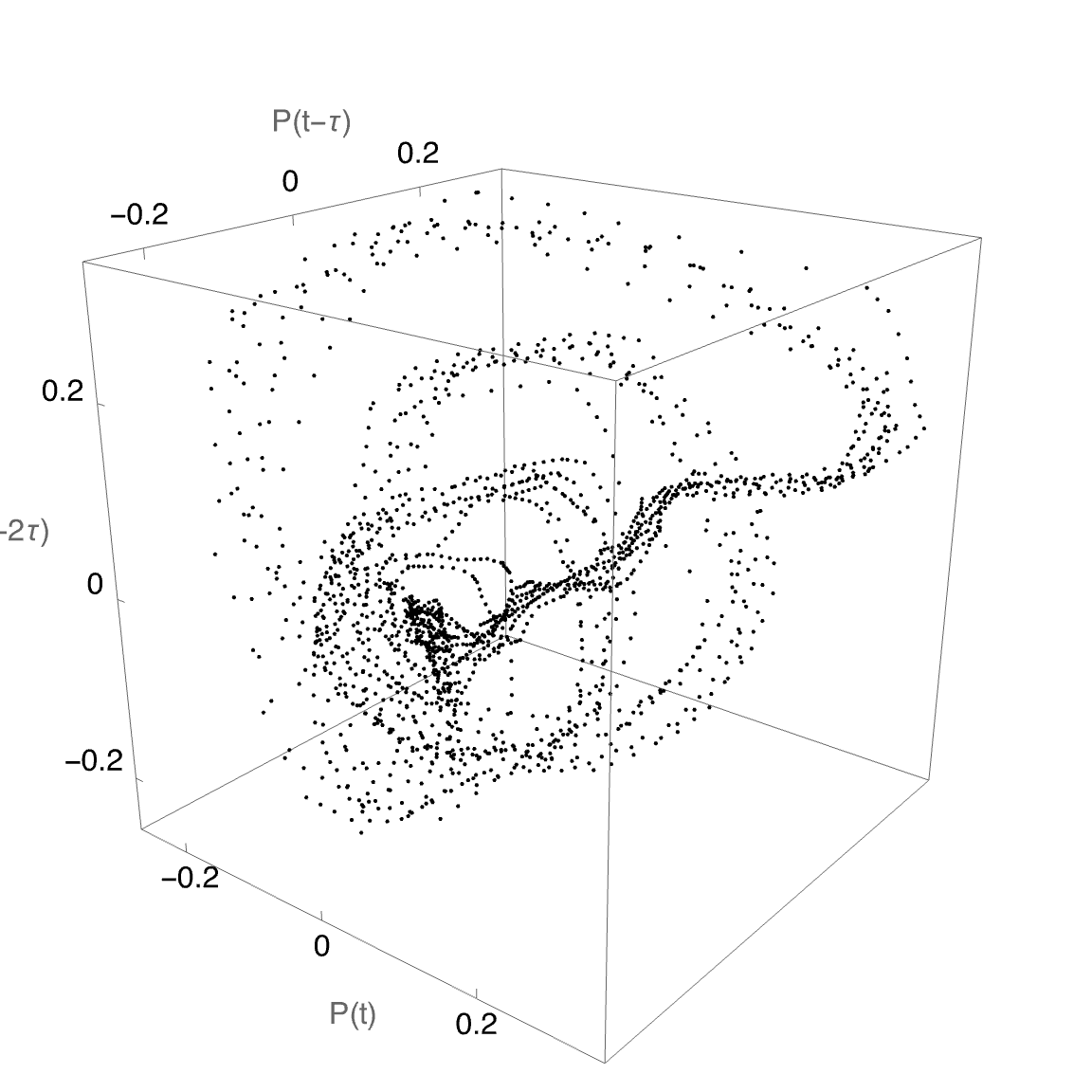}
\caption{The time series (a) of the pressure fluctuation values (in arbitrary units) as a function of their position in the data file. The three-dimensional time-delay embedding of the time series (b), with $\tau=5$.}
\label{data}
\end{figure}

This analysis reconsiders the dataset in Ref.~\cite{sciamarella1999topological}, namely, the pressure fluctuations recorded for the vowel /a/ as pronounced in Spanish. The time series and its time-delay embedding are displayed in Figs.\,\ref{data}(a,b). The file contains $1183$ points corresponding to a total interval of $t \approx 0.147~\text{s}$, with the measurements recorded at $8\;000$ samples per second. The BraMAH complex is reconstructed, so the number of 2-cells, their distribution, and labels differ from those in Ref.~\cite{sciamarella1999topological}. As expected, the homological features (three 1-holes) remain identical to those originally reported. This is a natural consequence of the formalism: the homology functor leads to invariants that are independent of the specific cell decomposition.

\begin{figure}[t]
(a)\includegraphics[width=0.8\columnwidth]{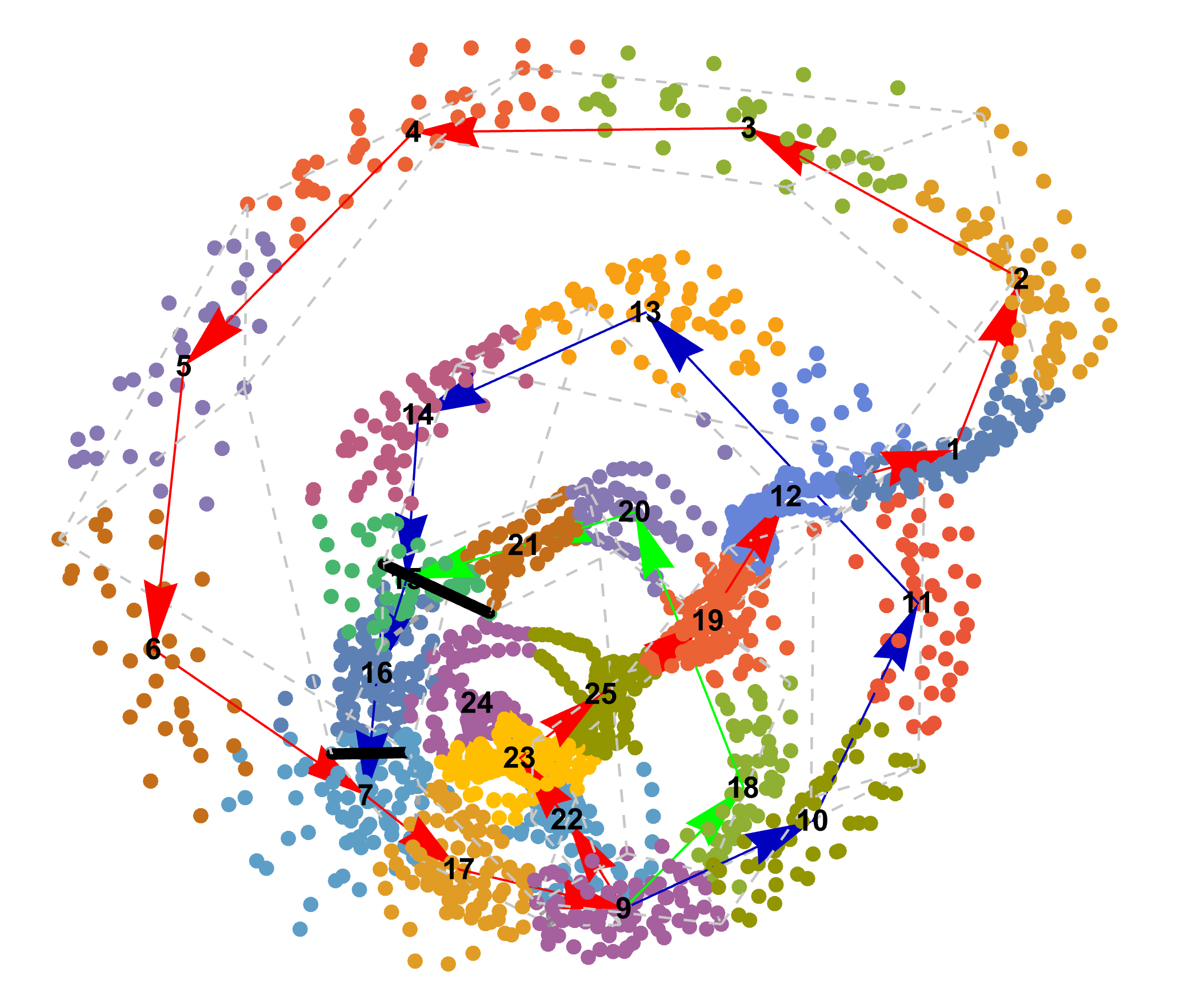}\\
(b)\includegraphics[width=0.8\columnwidth]{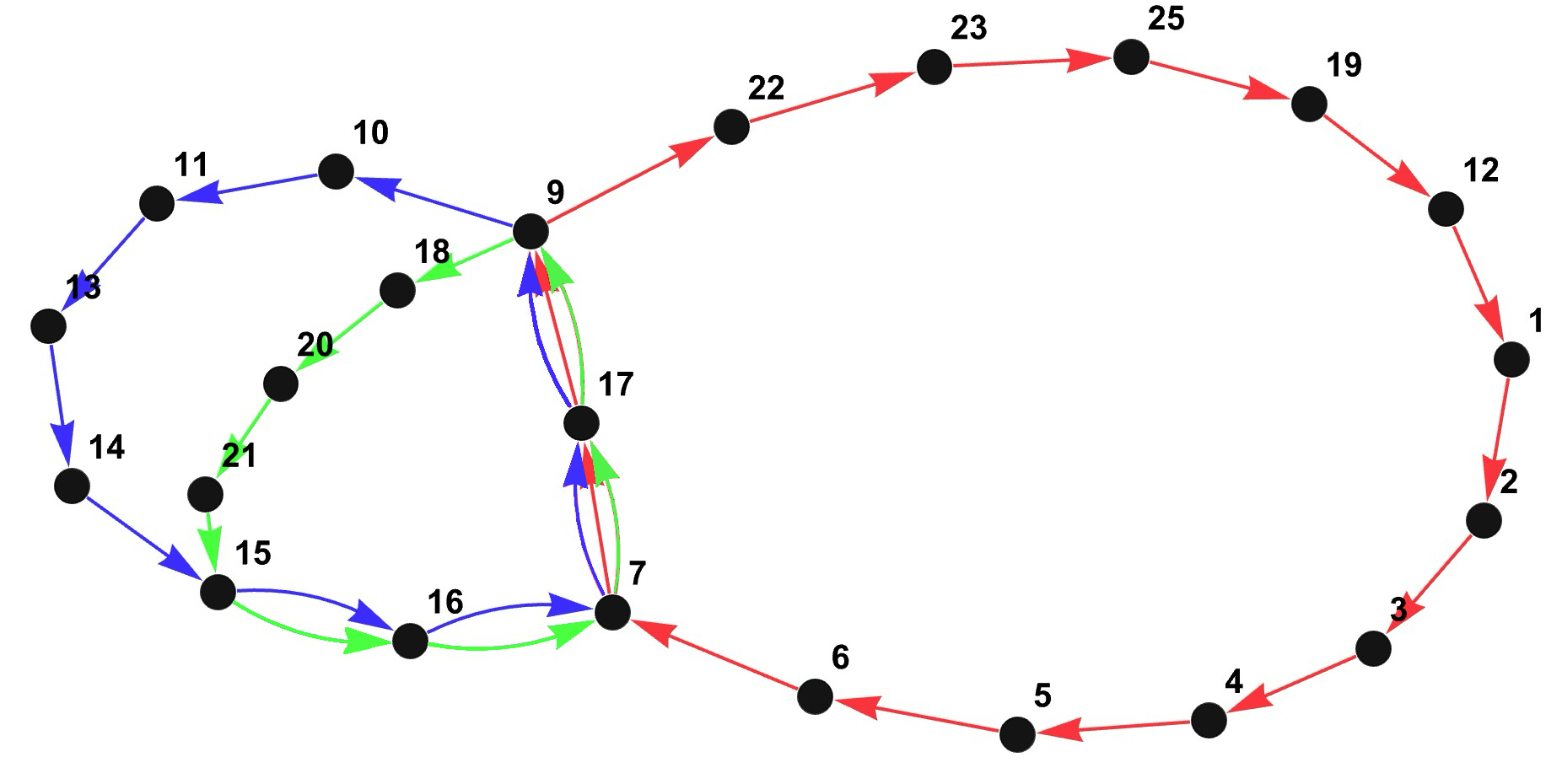}\\
(c)\includegraphics[width=0.8\columnwidth]{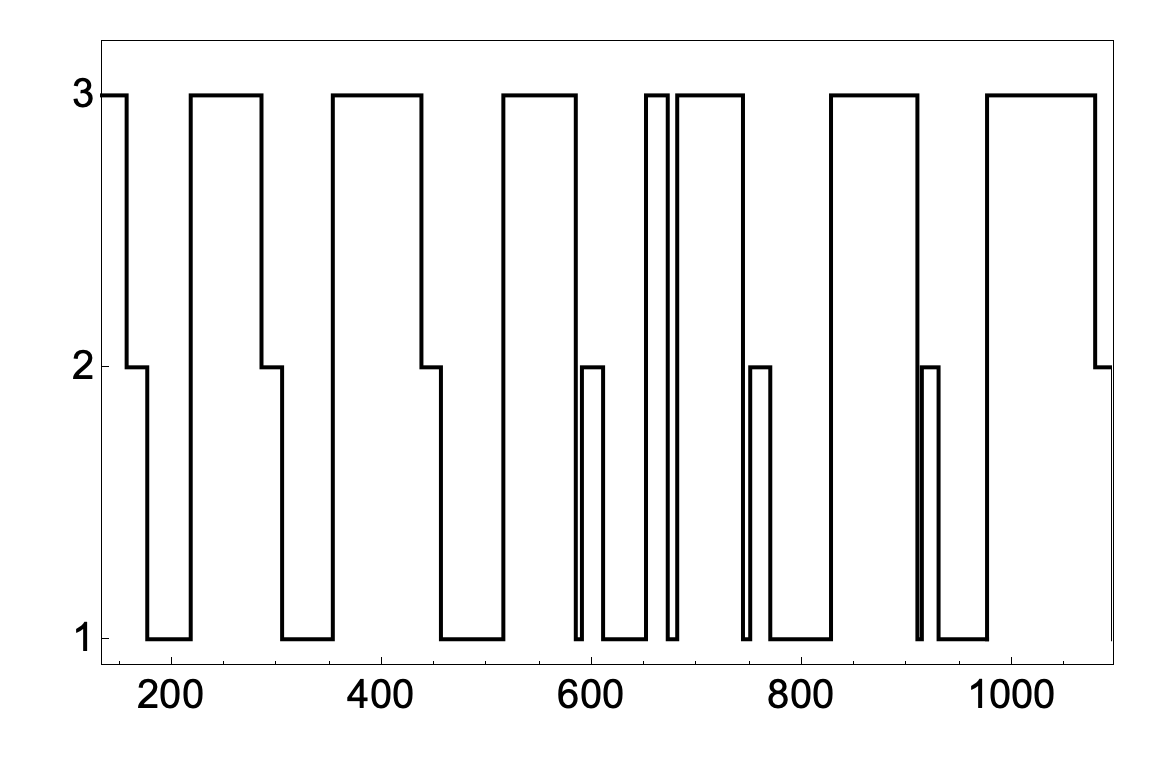}\\
\caption{(a) Point cloud segmented into different colors ($2$-cells) with $1$-cells shown as light-gray dashed lines. (b) Directed multigraph with the three generatexes in green, blue and red. (c) Generatex visitation sequence extracted from the time series shown in Fig.~\ref{data}a.
}
\label{speech}
\end{figure}

The analysis is conducted by applying the forgetful functor $U_{K}$, since homological properties were already discussed. The new aspect of the present work lies in the directed algebraic topology, which reveals three distinct generatex classes circling around two joining loci, indicated by thick black lines in Fig.~\ref{speech}(a):

\[
\begin{aligned}
G_1 &= \{7 \rightarrow 17 \rightarrow 9 \rightarrow 18 \rightarrow 20 \rightarrow 21 \rightarrow 15 \rightarrow 16 \rightarrow 7\}, \\
G_2 &= \{7 \rightarrow 17 \rightarrow 9 \rightarrow 10 \rightarrow 11 \rightarrow 13 \rightarrow 14 \rightarrow 15 \\
    &\quad \rightarrow 16 \rightarrow 7\}, \\
G_3 &= \{1 \rightarrow 2 \rightarrow 3 \rightarrow 4 \rightarrow 5 \rightarrow 6 \rightarrow 7 \rightarrow 17 \rightarrow 9 \\
    &\quad \rightarrow 22 \rightarrow 23 \rightarrow 25 \rightarrow 19 \rightarrow 12 \rightarrow 1\}.
\end{aligned}
\]
While $G_2$ is orientation-preserving, $G_1$ and $G_3$ are both orientation-reversing. The bonds, given by:
\[
\begin{aligned}
B_{123} &= G_1 \cap G_2 \cap G_3
        = \text{\{7 $\rightarrow$ 17 $\rightarrow$ 9\}},\\
B_{12}  &= G_1 \cap G_2
        = \text{\{15 $\rightarrow$ 16 $\rightarrow$ 7\}}.
\end{aligned}
\]
can be identified in the multigraph; see Fig.~\ref{speech}(b). Applying the Poincaré-edge operator $\mathcal{P}$ to each generatex, we obtain:
\[
\begin{aligned}
\mathcal{P}( G_1) &= ( \langle 21 \mid 15 \rangle, \langle 16 \mid 7 \rangle ),\\
\mathcal{P}( G_2) &= ( \langle 14 \mid 15 \rangle, \langle 16 \mid 7 \rangle ),\\
\mathcal{P}( G_3) &= ( \langle 6 \mid 7 \rangle ).
\end{aligned}
\]
While $\langle 14\mid 15\rangle$ and $\langle 21\mid 15\rangle$ traverse one joining locus, $\langle 16\mid 7\rangle$ and $\langle 6\mid 7\rangle$ traverse the other. The generatex classes are glued along a double bond, which becomes a triple one at node $7$. Each generatex class can be associated with a $1$-hole, yielding branches separated by voids. This contrasts with the Rössler attractor, where folding occurs without creating an internal void between the bonded generatex classes.

Fig.~\ref{speech}(c) shows the generatex visitation sequence. The signal alternates irregularly between the three generatexes $G_1,G_2,G_3$, with strongly variable residence times: approximately $320$ time units in $G_1$, $130$ in $G_2$, and $514$ in $G_3$. While transitions $G_3 \!\to\! G_1$ and self-concatenations $G_1 \!\to\! G_1$ dominate, the sequence $G_3 \!\to\! G_2 \!\to\! G_1$ appears only intermittently. This itinerant behavior across distinct generatex classes is the topological signature of the specific type of chaotic organization underlying the vocal-fold dynamics.

\subsection{Wind-driven double gyre}

\begin{figure*}[t]
\centering
(a)\includegraphics[width=0.9\columnwidth]{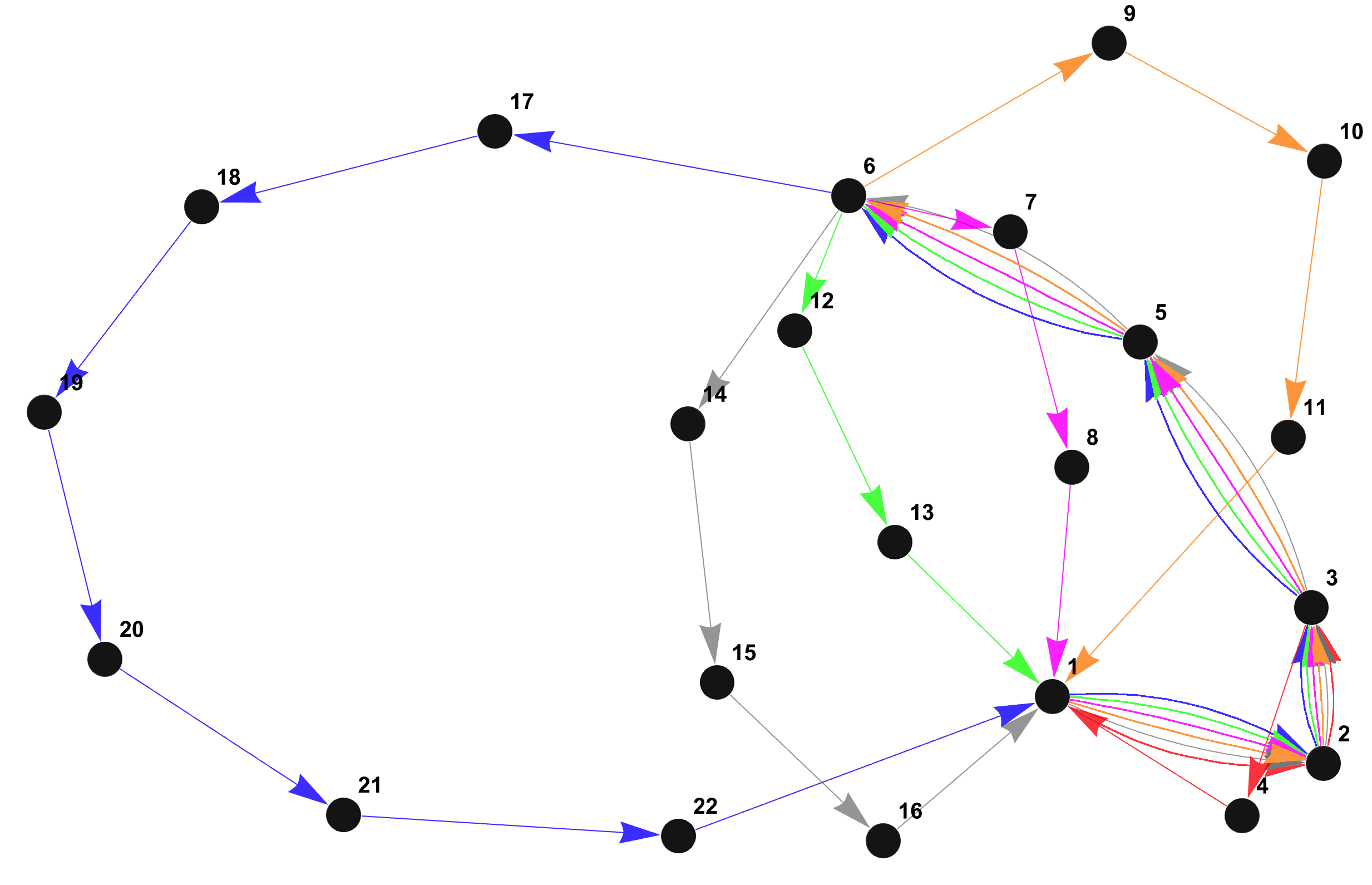}
(b)\includegraphics[width=0.85\columnwidth]{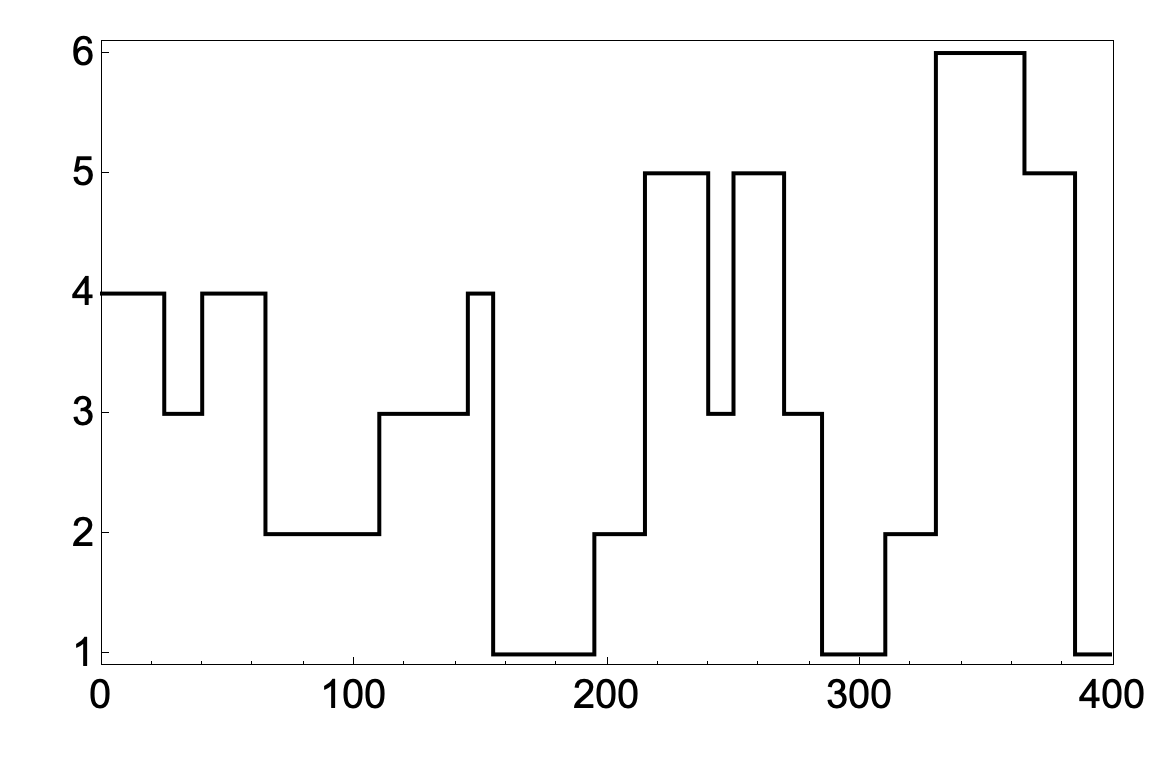}
\caption{
	(a) Directed multigraph obtained as the union of the six generatexes in red, blue, green, magenta, orange and gray for the templex of the aperiodically forced wind-driven double gyre example. 
	(b) Generatex visitation sequence extracted from the time series (in years) of an individual solution, shown in Fig.~15 of Ref.~\cite{charo2025topological}.
}
\label{fig:FMs}
\end{figure*}

Understanding the effects of time-dependent forcing on intrinsic ocean variability is crucial for modeling climate variability in general.  This section presents a numerical simulation reported in Ref.~\cite{charo2025topological} based on a model of intermediate complexity \citep{Ghil.Luc.2020} for the streamfunction of the wind-driven ocean circulation in midlatitudes \citep{Pier.Ghil.2021}. In the case of aperiodic wind-stress forcing, three junction loci (one joining locus and two splitting loci) are detected.  Note that top-cells are renumbered here as $\gamma_i \mapsto i$ ($i=1,\dots,6$) and $\sigma_j \mapsto 6+j$ ($j=1,\dots,16$). Six generatex classes are identified:
\[
\begin{aligned}
	\quad
	G_1 &= \{1 \rightarrow 2 \rightarrow 3 \rightarrow 4 \rightarrow 1\},\\
	G_2 &= \{1 \rightarrow 2 \rightarrow 3 \rightarrow 5 \rightarrow 6 \rightarrow 17 \rightarrow 18 \rightarrow 19\\
	&\qquad\rightarrow 20 \rightarrow 21 \rightarrow 22 \rightarrow 1\},\\
	G_3 &= \{1 \rightarrow 2 \rightarrow 3 \rightarrow 5 \rightarrow 6 \rightarrow 12 \rightarrow 13 \rightarrow 1\},\\
	G_4 &= \{1 \rightarrow 2 \rightarrow 3 \rightarrow 5 \rightarrow 6 \rightarrow 7 \rightarrow 8 \rightarrow 1\},\\
	G_5 &= \{1 \rightarrow 2 \rightarrow 3 \rightarrow 5 \rightarrow 6 \rightarrow 9 \rightarrow 10 \rightarrow 11 \rightarrow 1\},\\
	G_6 &= \{1 \rightarrow 2 \rightarrow 3 \rightarrow 5 \rightarrow 6 \rightarrow 14 \rightarrow 15 \rightarrow 16 \rightarrow 1\}.
\end{aligned}
\]
Only $G_2$ and $G_6$ are orientation-reversing. The directed multigraph is shown in Fig.~\ref{fig:FMs}(a). 
There are two bonds:
\[
\begin{aligned}
B_{123456} &= \bigcap_{i=1}^6 G_i 
= \text{\{1 $\rightarrow$ 2 $\rightarrow$ 3\}},\\
B_{23456} &= \bigcap_{i=2}^6 G_i= \text{\{1 $\rightarrow$ 2 $\rightarrow$ 3 $\rightarrow$ 5 $\rightarrow$ 6\}}.
\end{aligned}
\]
The valence of the bonds, given by the number of parallel edges, is
$6$ for $B_{123456}$ and $5$ for $B_{23456}$.  Applying the Poincaré-edge operator $\mathcal{P}$ to each generatex, we obtain:
\begin{alignat*}{2}
	\mathcal{P}(G_1) &= ( \langle 4 \mid 1 \rangle ), & \qquad \mathcal{P}(G_2) &= ( \langle 22 \mid 1 \rangle ), \\
	\mathcal{P}(G_3) &= ( \langle 13 \mid 1 \rangle ), & \qquad \mathcal{P}(G_4) &= ( \langle 8 \mid 1 \rangle ), \\
	\mathcal{P}(G_5) &= ( \langle 11 \mid 1 \rangle ), & \qquad \mathcal{P}(G_6) &= ( \langle 16 \mid 1 \rangle ).
\end{alignat*}
The joining locus being traversed is given by the 1-cell labeled $\langle 0,1 \rangle$ in Fig.~14(a) of Ref.~\cite{charo2025topological}. The generatex visitation label $\chi(t)$ defined in Section~\ref{sec:TMV} is shown in Fig.~\ref{fig:FMs}(b). The set of active TMVs grows over time: three generatex classes during the first 150 years, then two more (numbers 1 and 5), and the sixth towards the end of the inspected time series. Physically, the emergence (or disappearance) of generatex classes translates into the forcing-dependent accessibility of macroscopic flow configurations and provides a characterization of topological tipping phenomena in nonautonomous settings \cite{Pier.Ghil.2021,ashwin2026introduction}. The sustained vowel presents no tipping and in this sense is closer to the autonomous or periodically forced cases of the wind-driven double gyre discussed in Ref.~\cite{charo2025topological}, where the same type of chaotic behavior is maintained over the whole time window.

\section{Summary and conclusions}

The study of chaotic dynamics through topology is a program that can be formally broken into two parts: (i) the construction of the topological object from data; and (ii) the algebraic handling of this object. This article deals with (ii), assuming that (i) is adequately accomplished. The templex is placed within category theory, translating the problem of chaos topology into an algebraic one thanks to functors. Doing so has required providing an intrinsic definition of a BraMAH complex, an algebra of directed paths, an edge operator and its associated quotient, which is fundamentally different from the homology quotient. The result is a formal setting with two complementary levels of information carried by the functorial invariants and separable through forgetful functors. 

The undirected level (homology groups) describes the structural scaffold with its voids and orientability properties: these invariants allow one to distinguish, for instance, an underlying torus from other supporting shapes, regardless of the sequence in which regions are visited. If the scaffold corresponds to a pseudo-manifold, BraMAH detects and locates the junction loci. The directed level (generatex semigroups) captures the irreversible causal organization around joining loci, if they exist; these invariants are obtained through a directed algebra of paths, providing information about the specific type of chaotic dynamics at play. The classical mechanisms of chaotic phase-space organization naturally translate into this functorial framework: 

\begin{itemize}[leftmargin=\parindent]
\item[–] \textbf{Squeezing} is algebraically mapped to bonds which embody the topological overlap of distinct generatex classes; they are the regions where different futures are glued together, so that causal distinguishability is locally lost.
\item[–]  \textbf{Stretching} corresponds to the divergence of trajectories as they exit a bond and split into the distinct unbonded regions of their respective generatex classes, with the valence quantifying the number of non-equivalent directed paths available. 
\item[–]  \textbf{Folding} manifests as an orientation-reversing generatex class, captured by inspection of the orientability properties of each generatex studied as an independent and uniformly reoriented cell complex. Folding can occur with or without creating internal voids. 
\item[–]  \textbf{Tearing} actively splits the flow, typically segregating generatex classes and organizing them through different bonds, as observed in the Lorenz attractor. Tearing is absent in the Rössler case, which presents folding without tearing. 
\end{itemize}

One must beware that topological voids do not necessarily imply tearing; a quasi-periodic toroidal dynamics possesses two 1-holes and one 2-hole but exhibits no tearing mechanism whatsoever. The same toroidal support may host a regular dynamics, with a trivial generatex semigroup, or a chaotic one \cite{mosto2024templex}. Both functorial levels are thus relevant to unveil the nature of the topological space hosting the dynamics, and the type of directed topology inscribed upon it. The functorial touchstone for chaos is the existence of multiple, non-equivalent generatex classes. This is a non-metric criterion: the structural basis for sensitive dependence on initial conditions is algebraically encoded in the alternation between bonded and unbonded regions in phase space, inducing pushout diagrams in the category of semigroups.  The topological definition of chaos that results from this formalism is applicable to finite-time datasets, independently of whether they are obtained from numerical integrations of well-known chaotic systems, from experimental time series as in the voice production case, or from higher-dimensional models as in the wind-driven double gyre.

The speech example constitutes the first templex constructed from a fully experimental signal and the first application of the functorial formulation to voice dynamics. The resulting directed invariants provide an operational signature of chaos in voiced sound production that remained undetectable by homologies alone. Chaotic signals in voice production can arise from disparate physiological mechanisms \cite{Svec2025}. Because clinical practice frequently conflates structurally distinct nonlinear phenomena under generic perceptual terms like diplophonia or roughness, functorial invariants offer a rigorous tool for disambiguation. Future research applying TMV decompositions to empirical acoustic data from distinct clinical cohorts could reveal whether specific laryngeal pathologies impose unique topological signatures, ultimately enabling diagnosis through non-invasive acoustic recordings. On the other hand, the change from a voiced vowel to an unvoiced sound should translate into a global topological reconfiguration.

The climatic example shows how the formalism works with a four-dimensional phase space and nonautonomous settings. Topological tipping occurs at those instants when the set of available generatex classes changes. As discussed in Ref.~\cite{crisan2026introduction}, ensemble realizations of a chaotic system lead to the coexistence of several TMVs at a single instant, offering a topological counterpart to the changes in statistical moments driven by time-dependent forcing \citep{Maraldi.ea.2025}. In the context of the large-scale ocean circulation, these directed invariants capture the spontaneous reorganization of mass transport pathways, such as abrupt shifts in the separation of western boundary currents or transitions between distinct gyre regimes. Topological tipping thus encodes the switching of the ocean dynamics between structurally different, physically realizable states. Topological tipping points defined in the context of random attractors with a random templex require an extension of this formalism to noise-driven chaotic dynamics which is out of the scope of this work \cite{charo2021noise,charo2023random}.

Because the templex is constructed by processing any dynamical evolution as a point cloud, one can use the full set of invariants for the intercomparison of datasets, the cross-validation of models, and the benchmarking of simulated outputs against physical observations. In applications, one must work out a dictionary of the concrete meaning of the functorial invariants in terms of the specific problem, which may involve addressing broad concepts such as tipping points \cite{mosto2025templex,charo2025topological}, coherent structures \cite{bonel2025templex}, weather regimes \cite{strommen2023topological} or extreme events \cite{letellier2025extreme}. Functorial invariants can also help advance physics-informed machine learning. Traditional dimensionality reduction techniques, such as Principal Component Analysis or Empirical Orthogonal Decompositions, rely on projections that optimize variance but do not necessarily guarantee the preservation of the undirected and directed topology of the data. Similarly, while complex network approaches to time series analysis successfully capture sequential transitions, they often lack the formal topological scaffold required to define robust invariants. The architecture and cost functions of a data-driven model can be validated, or actively guided, based on their capacity to reproduce the generatex semigroup and bond structure of the original system or dataset. 

To conclude, this article presents a functorial formulation for chaos topology from data. Placing the templex within category theory discloses why homology must not be replaced by directed homology, what kind of object the templex is, and why two related but distinct functorial chains are required to translate finite-time dynamics into a complete set of topological invariants. This contribution brings together\,—\,while keeping them distinct\,—\,two notions that are traditionally difficult to reconcile: time and structure.\\

{\it{Acknowledgments}}
The author wishes to thank Michael Ghil and Pablo Amster for careful reading and insightful discussions, and Christophe Letellier for exchanges on templates. Funds from the ANR project TeMPlex ANR-23-CE56-0002 are acknowledged.

\appendix

\section{Homologies}
\label{app:homology}

A \emph{cell complex} is a topological space built by gluing cells of increasing dimension (points, segments, disks, and higher-dimensional analogues) along their boundaries~\cite{kinsey2012topology}. 
It provides a discrete, combinatorial representation of an underlying geometric object, and can be regarded as a topological scaffold: while the complex itself is not the object of interest, it supports algebraic constructions that reveal intrinsic properties of the space it represents. The dimension of a cell complex is defined as the highest dimension of its cells.

Given a finite cell complex $K$, the group of $k$-chains $C_k(K)$ is defined as the free abelian group generated by the $k$-cells of $K$,
\[
C_k(K)=\left\{\sum_i n_i\,\sigma_i^k \;\middle|\; n_i\in\mathbb{Z},\ \sigma_i^k \text{ a $k$-cell of } K\right\}, \]
together with a boundary operator
\[
\partial_k : C_k(K)\longrightarrow C_{k-1}(K),
\]
defined by assigning to each $k$-cell the formal sum of its $(k-1)$-dimensional faces, with integer coefficients determined by their relative orientations.  The graded family of chain groups together with the boundary operators
$\partial_k : C_k(K)\to C_{k-1}(K)$ that satisfy $\partial_{k-1}\circ\partial_k=0$
is denoted by $C_\bullet(K)=\{C_k(K)\}_{k\ge 0}$.
The group of $k$-cycles is the subgroup
\[
Z_k(K)=\ker\partial_k \subset C_k(K),
\]
while the group of $k$-boundaries is the subgroup
\[
B_k(K)=\mathrm{im}\,\partial_{k+1}\subset Z_k(K).
\]
The $k$-th homology group is then defined as the quotient
\[
H_k(K)=Z_k(K)/B_k(K).
\]
Intuitively, forming a quotient amounts to identifying objects that play the same structural role across different representations. A useful metaphor is that of musical notes: D3 and D4 are distinct pitches, separated by an octave, yet both belong to the same note class “D”, which might be known to some readers as “Re.” In homology, cycles that differ only by the attached higher–dimensional cells are identified in an analogous way, as they represent the same topological feature.

The generators of $H_k(K)$ correspond to intrinsic $k$-dimensional holes of the underlying space: they represent closed $k$-cycles that cannot be expressed as the boundary of any $(k+1)$-dimensional chain. In this sense, homology detects topological features that are invariant under changes of the cell decomposition, capturing global properties such as connected components ($0$-holes), tunnels ($1$-holes), and higher-dimensional cavities ($k\ge 2$). A generator of $H_k(K)$ can be represented by any closed $k$-cycle within its equivalence class.

Betti numbers $\beta_k$, defined as the rank of the homology groups $H_k$, count the number of independent $k$-dimensional holes. However, these homological ranks alone are insufficient to fully characterize the topology of the support; for instance, a cylinder and a Möbius band share identical Betti numbers. To detect the orientability properties of the complex, BraMAH evaluates the \emph{orientability chain}, defined as the sum of the boundaries of all top-dimensional cells~\cite{sciamarella2001unveiling}. As illustrated in Table~\ref{tab:betti_surfaces}, the evaluation of these chains provides the missing topological distinction, allowing one to discriminate between orientable and non-orientable manifolds.

\begin{table}[ht]
  \centering
  \caption{Betti numbers $\beta_k$ ($k=0,1,2$) and orientability properties for classic surfaces. 
  }
  \label{tab:betti_surfaces}
  \begin{tabular}{lcccc} 
    \hline \hline
    \\[-0.3cm]
    Cell complex & $\beta_0$ & $\beta_1$ & $\beta_2$ & Orientability \\[0.1cm]
    \hline 
    \\[-0.3cm]
    Sphere ($S^2$) & $1$ & $0$ & $1$ & Orientable \\ 
    Cylinder & $1$ & $1$ & $0$ & Orientable \\ 
    Möbius band & $1$ & $1$ & $0$ & Non-orientable \\ 
    Torus ($T^2$) & $1$ & $2$ & $1$ & Orientable \\ 
    Klein bottle & $1$ & $1$ & $0$ & Non-orientable \\[0.1cm]
    \hline \hline
  \end{tabular}
\end{table}

\section{The Rössler and Lorenz attractors}
\label{app:RL}

\begin{figure}[!htbp]
    \centering
    \begin{subfigure}[b]{0.47\columnwidth}
        \centering
        \includegraphics[width=\linewidth]{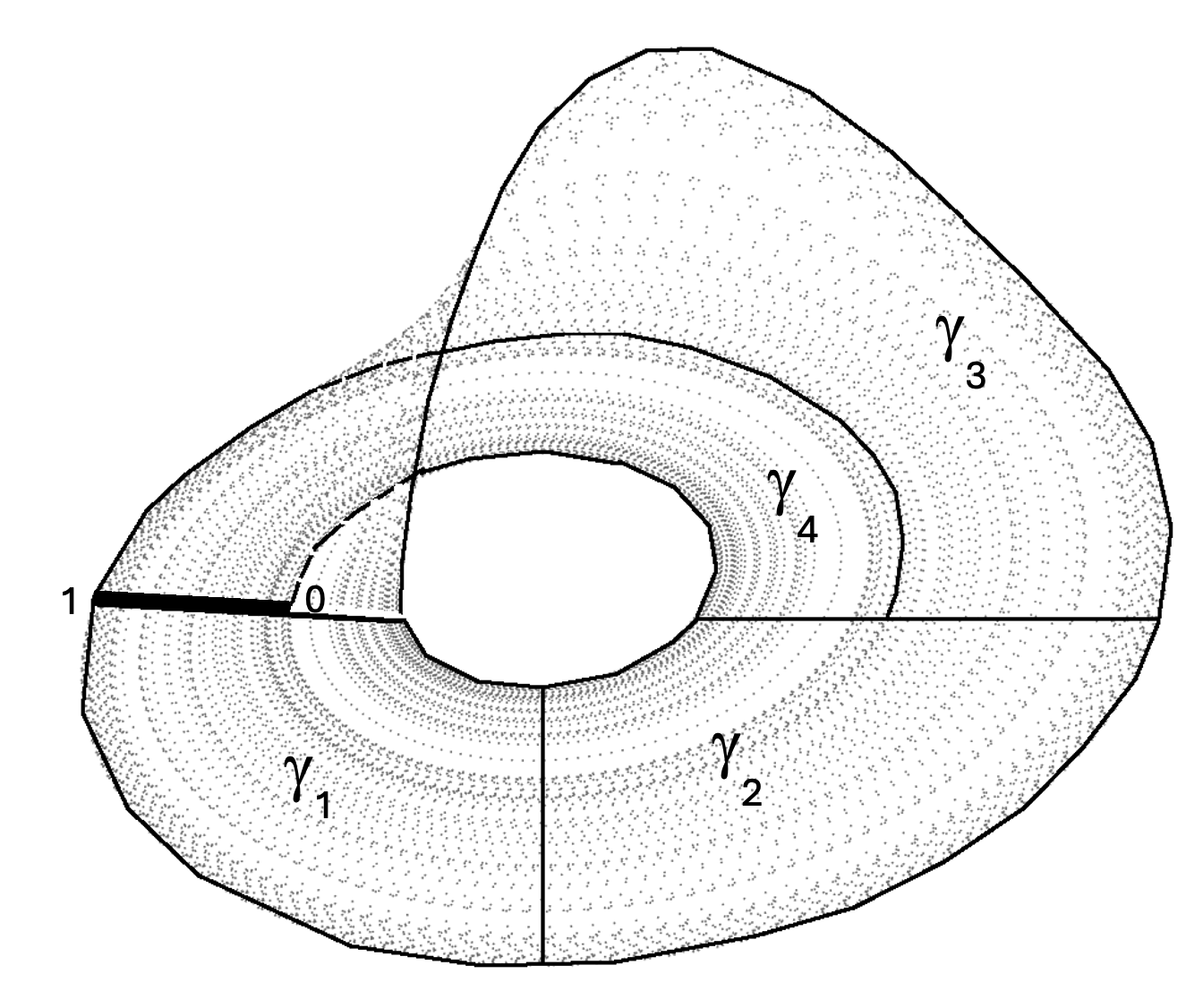}
        \caption{}
        \label{RosBraMAH}
    \end{subfigure}
    \hfill
    \begin{subfigure}[b]{0.47\columnwidth}
        \centering
        \includegraphics[width=\linewidth]{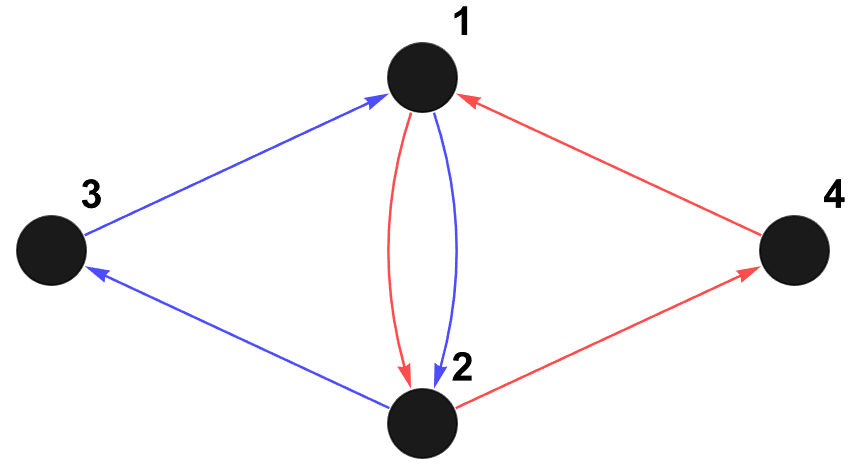}
        \caption{}
        \label{fig:rossler_multi}
    \end{subfigure}
    \\[0.8em]
    \begin{subfigure}[b]{0.75\columnwidth}
        \centering
        \includegraphics[width=\linewidth]{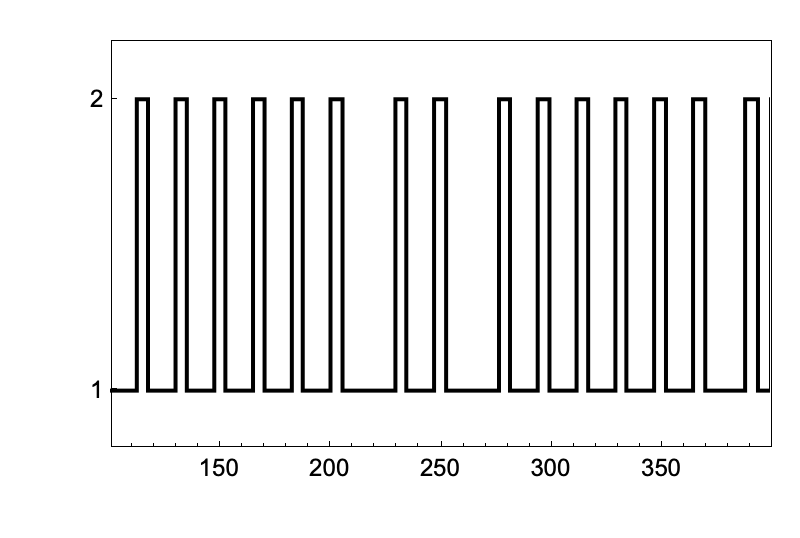}
        \caption{}
        \label{fig:rossler_visits}
    \end{subfigure}
\caption{The Rössler attractor. (a) BraMAH complex $K_R$; the junction locus is indicated by a thick line. (b) Directed multigraph for the Rössler templex $\mathcal{T}_R$. (c) Generatex visitation sequence.}
    \label{fig:rossler}
\end{figure}

\begin{figure}[!htbp]
    \centering
    \begin{subfigure}[b]{0.47\columnwidth}
        \centering
        \includegraphics[width=\linewidth]{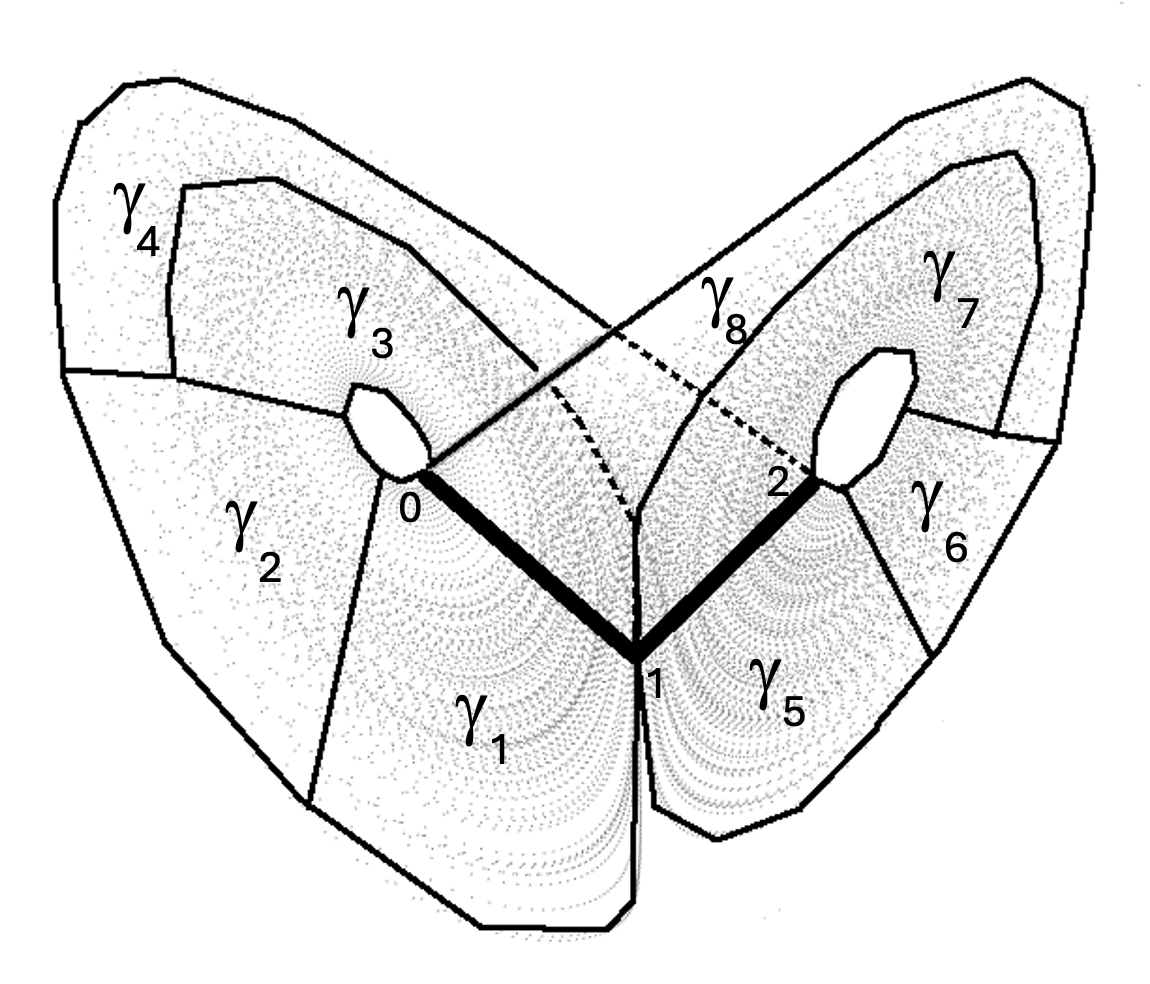}
        \caption{}
        \label{LorBraMAH}
    \end{subfigure}
    \hfill
    \begin{subfigure}[b]{0.47\columnwidth}
        \centering
        \includegraphics[width=\linewidth]{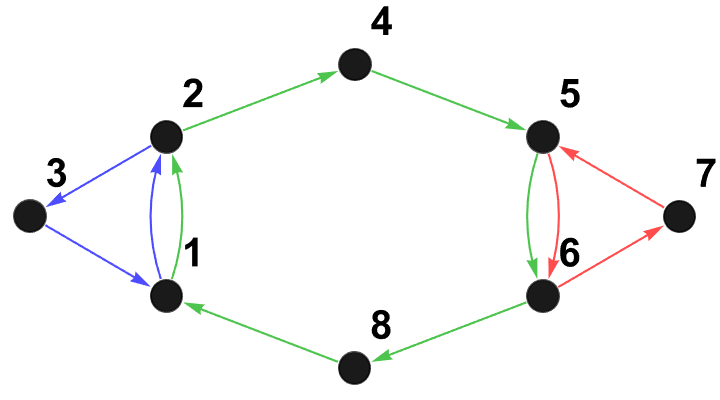}
        \caption{}
        \label{fig:lorenz_multi}
    \end{subfigure}
    \\[0.8em]
    \begin{subfigure}[b]{0.75\columnwidth}
        \centering
        \includegraphics[width=\linewidth]{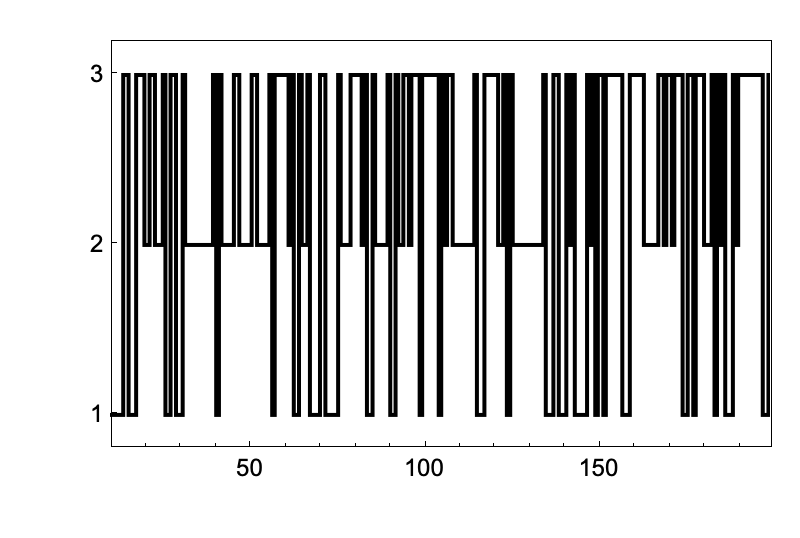}
        \caption{}
        \label{fig:lorenz_visits}
    \end{subfigure}
\caption{The Lorenz attractor. (a) BraMAH complex $K_L$; the two joining loci are indicated by thick lines. (b) Directed multigraph for the Lorenz templex $\mathcal{T}_L$. (c) Generatex visitation sequence.
}
    \label{fig:lorenz}
\end{figure}

This appendix presents the functorial analysis of two paradigmatic attractors, namely, the Rössler and the Lorenz attractors. Figures~\ref{fig:rossler} and \ref{fig:lorenz} display, for each of them, the BraMAH complexes $K_R$ and $K_L$, the directed multigraphs and generatex visitation sequences. To see how the same invariants are obtained from templexes with different cell decompositions, please see Refs.~\cite{charo2022templex,sciamarella2024new}.

\paragraph{Rössler attractor}

The BraMAH complex $K_R$ is built from four $2$-cells, $\gamma_1,\dots,\gamma_4$, and its junction locus consists of the single $1$-cell $\langle 0,1\rangle$ (Fig.~\ref{RosBraMAH}), where $\langle i,j\rangle$ denotes the $1$-cell whose endpoint $0$-cells are $\langle i \rangle$ and $\langle j \rangle$. Its homology groups are
\[
H_0(K_R)\cong\mathbb{Z},\quad H_1(K_R)\cong\mathbb{Z},\quad
H_k(K_R)=0\;\text{for }k\ge 2,
\]
the single generator of $H_1(K_R)$ representing the central void around which trajectories wind.

Adding the directed graph $\mathcal{G}_R$, we get $\mathcal{T}_R$. The digraph consists of nodes $N_R = \{1, 2, 3, 4\}$ corresponding to the 2-cells  $\{\gamma_1, \gamma_2, \gamma_3, \gamma_4\}$ and directed edges $E_R = \{(1,2), (2,3), (2,4), (3,1), (4,1)\}$.  Node 1 ($\gamma_1$) acts as an outgoing node (2-cell), receiving flow from ingoing nodes 3 and 4.  Two generatex classes are identified:
\[
\begin{aligned}
G_1 &= \{1 \rightarrow 2 \rightarrow 3 \rightarrow 1\}, \\
G_2 &= \{1 \rightarrow 2 \rightarrow 4 \rightarrow 1\}.
\end{aligned}
\] 
When considered as isolated surfaces, $G_1$ is homeomorphic to a standard cylinder ($G_1$ is orientation-preserving), while $G_2$ forms a Möbius band ($G_2$ is orientation-reversing). Their union is represented with the directed multigraph (Fig.~\ref{fig:rossler_multi}), where the single bond $B_{12} = G_1 \cap G_2 = \{1 \rightarrow 2\}$ has a valence of 2, evidenced by the two parallel edges connecting nodes 1 and 2. 

Applying the Poincaré-edge operator $\mathcal{P}$ to each generatex, we obtain:
\[
\begin{aligned}
\mathcal{P}(G_1) &= ( \langle 3 \mid 1 \rangle ),\\
\mathcal{P}(G_2) &= ( \langle 4 \mid 1 \rangle ).
\end{aligned}
\]
Both Poincaré edges traverse the single joining locus, the 1-cell $\langle 0,1\rangle$ of $K_R$.

Applying the forgetful functor $U_{\mathcal{G}}$ to $\mathcal{T}_R$ recovers the classical chain complex $(C_\bullet(K_R), \partial)$, with $H_1(K_R) \cong \mathbb{Z}$ encoding the single central void. Applying $U_K$ yields the directed path algebra $(C_{\mathcal{G}}(\mathcal{T}_R), \mathcal{P})$, and the generatex functor $F_{\mathrm{Gen}}$ maps it to the semigroup with two generatex classes.  The outer, orientation-reversing generatex class ($G_2$) is associated with Rössler's folding mechanism. 

It is worth noting that directed homology (dihomology), conceived to detect directed loops around topological voids, would not detect this mechanism. From a dihomological perspective, the Rössler branched manifold and a simple directed annulus (such as those arising from certain quasiperiodic flows) are equivalent, both enclosing a single 1-hole. The same applies to recurrence-based topological methods (e.g., cycling signatures from persistent homology \cite{bauer2023cycling}), which cannot separate the folding from the non-folding region.

The TMV decomposition allows computing the visitation sequence, which exhibits irregular alternation between the two generatex labels, with variable residence times $|I_k|$; see Fig.~\ref{fig:rossler_visits}. The two classes have an almost equal number of visits, but markedly unequal residence times: the total time spent in $G_1$ is approximately $220$ time units, with a mean visit duration of $14$, while $G_2$ accumulates only $79$ time units with a mean of $5$. The asymmetry in residence times reflects the well-known geometry of the spiral Rössler attractor: trajectories wind slowly outward through many revolutions before reaching the fold, i.e., long visits in $G_1$, and are then rapidly reinjected through the outer branch back to the inner spiral, i.e., short visits in $G_2$. For this realization, the mean visit duration ratio $\langle |I|_{G_1} \rangle / \langle |I|_{G_2} \rangle \simeq 2.8$ captures the separation between the slow spiraling and the fast folding mechanism.

\paragraph{Lorenz attractor}

The butterfly-shaped structure in $K_L$ (Fig.~\ref{LorBraMAH}) presents two joining loci at the $1$-cells $\langle 0,1\rangle$ and $\langle 1,2\rangle$. There are two $1$-holes, one per wing, so that:
\[
H_0(K_L)\cong\mathbb{Z},\quad H_1(K_L)\cong\mathbb{Z}^2,\quad
H_k(K_L)=0\;\text{for }k\ge 2.
\]

In $\mathcal{T}_L=(K_L,\mathcal{G}_L)$, the digraph $\mathcal{G}_L$ has $N_L = \{1, \dots, 8\}$ from the eight 2-cells $\{\gamma_1, \dots, \gamma_8\}$, and directed edges $E_L = \{(1,2), (2,3), (2,4), (3,1),$ $(4,5), (5,6), (6,7),$ $(7,5), (6,8), (8,1)\}$. There are two outgoing nodes 1 and 5, and four ingoing ones (3, 4, 7 and 8). Three orientation-preserving generatexes are identified:
\[
\begin{aligned}
G_1 &= \{1 \rightarrow 2 \rightarrow 3 \rightarrow 1\}, \\
G_2 &= \{5 \rightarrow 6 \rightarrow 7 \rightarrow 5\}, \\
G_3 &= \{1 \rightarrow 2 \rightarrow 4 \rightarrow 5 \rightarrow 6 \rightarrow 8 \rightarrow 1 \}.
\end{aligned}
\] 
$G_1$ and $G_2$ are of order $1$, each forming a stripex, while $G_3$ is of order $2$, forming two stripexes, yielding a total of four stripexes. The  union of the three generatexes forms a directed multigraph (Fig.~\ref{fig:lorenz_multi}), where the two bonds
\[
\begin{aligned}
B_{13} &= G_1 \cap G_3 = \{1 \rightarrow 2\},\\
B_{23} &= G_2 \cap G_3 = \{5 \rightarrow 6\},
\end{aligned}
\]
have a valence of 2, while $G_1$ and $G_2$ share no bond.  Applying the Poincaré-edge operator $\mathcal{P}$ to each generatex, we obtain:
\[
\begin{aligned}
\mathcal{P}(G_1) &= ( \langle 3 \mid 1 \rangle ),\\
\mathcal{P}(G_2) &= ( \langle 7 \mid 5 \rangle ),\\
\mathcal{P}(G_3) &= ( \langle 4 \mid 5 \rangle, \langle 8 \mid 1 \rangle ).
\end{aligned}
\]
The Poincaré edges traverse the two joining loci: $\langle 3 \mid 1 \rangle$ and $\langle 8 \mid 1 \rangle$ go through the 1-cell $\langle 0,1\rangle$, while $\langle 4 \mid 5 \rangle$ and $\langle 7 \mid 5 \rangle$ through $\langle 1,2\rangle$.

The forgetful functor $U_{\mathcal{G}}$ recovers the two central voids of the two wings. $F_{\mathrm{Gen}}$ yields the semigroup with the three generatex classes. The partial bond structure captures the tearing mechanism around the saddle fixed point and gives two independent pushout diagrams,
\[
\begin{tikzcd}[column sep=large, row sep=large]
B_{13} \arrow[r, hook] \arrow[d, hook']
  & G_1 \arrow[d] \\
G_3 \arrow[r]
  & G_1 \sqcup_{B_{13}} G_3
\end{tikzcd}
\qquad
\begin{tikzcd}[column sep=large, row sep=large]
B_{23} \arrow[r, hook] \arrow[d, hook']
  & G_2 \arrow[d] \\
G_3 \arrow[r]
  & G_2 \sqcup_{B_{23}} G_3
\end{tikzcd}
\]
\noindent both in $\mathbf{Sem}$. A direct transition from $G_1$ to $G_2$ (or vice versa) is thus forbidden. 

The generatex visitation sequence is shown in Fig.~\ref{fig:lorenz_visits}. $G_1$ and $G_2$ each record $26$ visits, while $G_3$ records $51 = 26 + 25$ visits, accumulating the largest total time ($82.8$ time units). For this specific realization, the mean residence times per visit are $\langle |I|_{G_1} \rangle = 1.7$, $\langle |I|_{G_2} \rangle = 2.4$, and $\langle |I|_{G_3} \rangle = 1.6$. The disparity between the two order-$1$ classes exemplifies how a symmetry in the governing equations does not preclude asymmetries in the residence times of finite-time realizations evaluated over a topological partition. It is interesting that analogous asymmetries in residence times between alternating regimes of a system exist, for instance, in climate and weather variability: the well-known El Niño-Southern Oscillation phenomenon in the tropical Pacific has longer and fewer El Niño vs. La Niña episodes measured in years \citep{Philander1990}. The blocking vs. zonal flow variability in the Northern Hemisphere midlatitudes has the same asymmetry but measured in weeks \citep{Ghil.Luc.2020}.

\bibliography{CatM}

@article{whitehead1949combinatorial,
  title={Combinatorial homotopy. I},
  author={Whitehead, John Henry Constantine},
  journal={Bulletin of the American Mathematical Society},
  volume={55},
  number={3},
  pages={213--245},
  year={1949},
  publisher={American Mathematical Society}
}

@article{karniadakis2021physics,
  title={Physics-informed machine learning},
  author={Karniadakis, George Em and Kevrekidis, Ioannis G and Lu, Lu and Perdikaris, Paris and Wang, Sifan and Yang, Liu},
  journal={Nature Reviews Physics},
  volume={3},
  number={6},
  pages={422--440},
  year={2021},
  publisher={Nature Publishing Group UK London}
}

@article{ashwin2026introduction,
  title={Introduction to Focus Issue: Nonautonomous dynamical systems: Theory, methods, and applications},
  author={Ashwin, Peter and Feudel, Ulrike and Ghil, Michael and Lehnertz, Klaus and Ortega, Juan-Pablo and Rasmussen, Martin},
  journal={Chaos: An Interdisciplinary Journal of Nonlinear Science},
  volume={36},
  number={4},
  year={2026},
  publisher={AIP Publishing}
}

@article{crisan2026introduction,
  author = {Crisan, Dan and Galatolo, Stefano and Ghil, Michael and Pierini, Stefano and Sciamarella, Denisse and Tél, Tamás},
  title = {Introduction to the Focus Issue: Nonautonomous dynamics in the climate sciences},
  journal = {Chaos: An Interdisciplinary Journal of Nonlinear Science},
  volume = {36},
  number = {040403},
  pages = {040403},
  year = {2026},
  doi = {10.1063/5.0324361},
  url = {https://doi.org/10.1063/5.0324361}
}

@article{mosto2024templex,
  title={Templex-based dynamical units for a taxonomy of chaos},
  author={Mosto, Caterina and Char{\'o}, Gisela D and Letellier, Christophe and Sciamarella, Denisse},
  journal={Chaos: An Interdisciplinary Journal of Nonlinear Science},
  volume={34},
  number={11},
  year={2024},
  publisher={AIP Publishing}
}

@article{brunton2020machine,
  title={Machine learning for fluid mechanics},
  author={Brunton, Steven L and Noack, Bernd R and Koumoutsakos, Petros},
  journal={Annual review of fluid mechanics},
  volume={52},
  number={1},
  pages={477--508},
  year={2020},
  publisher={Annual Reviews}
}

@article{Let10b,
  author  = {Letellier, Christophe and Messager, V{\'e}ronique},
  title   = {Influences on {O}tto {E}. {R}{\"o}ssler's earliest paper on chaos},
  journal = {International Journal of Bifurcation and Chaos},
  year    = {2010},
  volume  = {20},
  number  = {11},
  pages   = {3585--3616},
  doi     = {10.1142/S021812741002797X}
}

@book{MacLane1971,
  author    = {Saunders {Mac Lane}},
  title     = {Categories for the Working Mathematician},
  series    = {Graduate Texts in Mathematics},
  volume    = {5},
  publisher = {Springer},
  year      = {1971},
  isbn      = {978-0387900355},
  doi       = {10.1007/978-1-4757-4721-8}
}

@article{EilenbergMacLane1945,
  author    = {Samuel Eilenberg and Saunders {Mac Lane}},
  title     = {General theory of natural equivalences},
  journal   = {Transactions of the American Mathematical Society},
  volume    = {58},
  number    = {2},
  pages     = {231--294},
  year      = {1945},
  doi       = {10.2307/1990284}
}

@article{mosto2025templex,
  title={A templex-based study of the Atlantic Meridional Overturning Circulation dynamics in idealized chaotic models},
  author={Mosto, Caterina and Char{\'o}, Gisela D and S{\'e}vellec, Florian and Tandeo, Pierre and Ruiz, Juan J and Sciamarella, Denisse},
  journal={Chaos: An Interdisciplinary Journal of Nonlinear Science},
  volume={35},
  number={1},
  year={2025},
  publisher={AIP Publishing}
}

@Article{charo2022templex,
  author    = {Char{\'o}, Gisela D. and Letellier, Christophe and Sciamarella, Denisse},
  journal   = {Chaos: An Interdisciplinary Journal of Nonlinear Science},
  title     = {Templex: A bridge between homologies and templates for chaotic attractors},
  year      = {2022},
  number    = {8},
  pages     = {083108},
  volume    = {32},
  publisher = {AIP Publishing LLC},
}

@article{charo2021noise,
  title={Noise-driven topological changes in chaotic dynamics},
  author={Char{\'o}, Gisela D and Chekroun, Micka{\"e}l D and Sciamarella, Denisse and Ghil, Michael},
  journal={Chaos: An Interdisciplinary Journal of Nonlinear Science},
  volume={31},
  number={10},
  pages={103115},
  year={2021},
  publisher={AIP Publishing LLC}
}

@article{charo2021topological,
title={Topological colouring of fluid particles unravels finite-time coherent sets}, 
author={Charó, Gisela D. and Artana, Guillermo and Sciamarella, Denisse},
volume={923}, DOI={10.1017/jfm.2021.561}, 
journal={J. Fluid Mech.},
publisher={Cambridge University Press},
year={2021},
pages={A17}
}

@article{mangiarotti2014modelisation,
  title={Mod{\'e}lisation globale et caract{\'e}risation topologique de dynamiques environnementales: de l’analyse des enveloppes fluides et du couvert de surface de la Terre {\`a} la caract{\'e}risation topolodynamique du chaos},
  author={Mangiarotti, Sylvain and others},
  journal={Habilitation to direct Researches, Universit{\'e} de Toulouse},
  volume={3},
  year={2014}
}

@article{lefranc2006alternative,
  title={Alternative determinism principle for topological analysis of chaos},
  author={Lefranc, Marc},
  journal={Physical Review E—Statistical, Nonlinear, and Soft Matter Physics},
  volume={74},
  number={3},
  pages={035202},
  year={2006},
  publisher={APS}
}

@Book{gilmore2012topology,
  author    = {Gilmore, Robert and Lefranc, Marc},
  publisher = {John Wiley \& Sons},
  title     = {{The Topology of Chaos: Alice in Stretch and Squeezeland}},
  year      = {2012},
}

@ARTICLE{gilmore1998topological,
title={Topological analysis of chaotic dynamical systems},
   author       = "R. Gilmore",
   year         = "1998",
   journal      = "Rev.\ Mod.\ Phys.",
   volume       = "4",
   pages        = "1455",
}

@article{edelsbrunnerHarer2010,
  title = {Computational Topology: An Introduction},
  author = {Edelsbrunner, Herbert and Harer, John L.},
  journal = {American Mathematical Society},
  year = {2010},
  note = {See Chapter~III for Čech and Vietoris–Rips complexes.}
}

@ARTICLE{sciamarella1999topological,
title={Topological structure of chaotic flows from human speech data},
   author       = "D. Sciamarella  and G. B. Mindlin",
   year         = "1999",
   journal      = "Phys.\ Rev.\ Lett.",
   volume       = "82",
   pages        = "1450",
}

@ARTICLE{sciamarella2001unveiling,
 title={Unveiling the topological structure of chaotic flows from data},
   author       = "D. Sciamarella  and G. B. Mindlin",
   year         = "2001",
   journal      = "Phys.\ Rev.\ E",
   volume       = "64",
   pages        = "036209",
}

@ARTICLE{lorenz1963lorenz,
title="Deterministic nonperiodic flow",
   author       = "E. N. Lorenz",
   year         = "1963",
   journal      = "J.\ Atmos.\ Sci.",
   volume       = "20",
   pages        = "130--141",
}

@ARTICLE{charo2020topology,
title={Topology of dynamical reconstructions from {L}agrangian data},
   author       = "G. D. Char\'o and G. Artana and D. Sciamarella",
   year         = "2020",
   journal      = "Phys.\ D",
   volume       = "405",
   pages        = "132371",
}

@article{letellier2025extreme,
  title={To be an extreme event or not: That is the question},
  author={Letellier, Christophe and Kamdjeu Kengne, L{\'e}andre and Zhao, Manyu and Minati, Ludovico},
  journal={Chaos: An Interdisciplinary Journal of Nonlinear Science},
  volume={35},
  number={7},
  year={2025},
  publisher={AIP Publishing}
}

@article{bonel2025templex,
  title={Templex for Lagrangian dynamics in the Southwestern Atlantic},
  author={Bonel, Juan Cruz and Bodnariuk, Nicol{\'a}s and Char{\'o}, Gisela D and Letellier, Christophe and Guinet, Christophe and Saraceno, Mart{\'\i}n and Sciamarella, Denisse},
  journal={Chaos: An Interdisciplinary Journal of Nonlinear Science},
  volume={35},
  number={10},
  year={2025},
  publisher={AIP Publishing}
}

@article{charo2025topological,
  title={Topological variability modes of the wind-driven ocean circulation},
  author={Char{\'o}, Gisela D and Sciamarella, Denisse and Ruiz, Juan and Pierini, Stefano and Ghil, Michael},
  journal={arXiv e-prints},
  pages={arXiv--2502},
  year={2025}
}

@article{strommen2023topological,
  title={A topological perspective on weather regimes},
  author={Strommen, Kristian and Chantry, Matthew and Dorrington, Joshua and Otter, Nina},
  journal={Climate Dynamics},
  volume={60},
  number={5},
  pages={1415--1445},
  year={2023},
  publisher={Springer}
}

@Article{Ghil.Luc.2020,
  author    = {Michael Ghil and Valerio Lucarini},
  journal   = {Reviews of Modern Physics},
  title     = {The physics of climate variability and climate change},
  year      = {2020},
  number    = {3},
  pages     = {035002},
  volume    = {92},
  doi       = {10.1103/revmodphys.92.035002},
  publisher = {American Physical Society ({APS})},
}

@Article{Pier.Ghil.2021,
  author    = {Pierini, S. and Ghil, M.},
  journal   = {Sci. Rep.},
  title     = {Tipping points induced by parameter drift in an excitable ocean model},
  year      = {2021},
  number    = {1},
  pages     = {1--14},
  volume    = {11},
  publisher = {Nature Publishing Group},
  url       = {https://rdcu.be/clp5V},
}

@Book{kinsey2012topology,
  author    = {Kinsey, L. Christine},
  publisher = {Springer Science \& Business Media},
  title     = {{Topology of Surfaces}},
  year      = {2012},
}

@article{Ros76c,
author = {O. E. R{\"o}ssler},
title = "An equation for continuous chaos",
journal = "Physics Letters A",
volume = "57",
number = "5",
pages = "397-398",
year = "1976",
doi = {10.1016/0375-9601(76)90101-8},
}

@article{ghil2023dynamical,
  title={Dynamical systems, algebraic topology and the climate sciences},
  author={Ghil, Michael and Sciamarella, Denisse},
  journal={Nonlinear Processes in Geophysics},
  volume={30},
  number={4},
  pages={399--434},
  year={2023},
  publisher={Copernicus Publications G{\"o}ttingen, Germany}
}

@article{charo2023random,
  title={Random templex encodes topological tipping points in noise-driven chaotic dynamics},
  author={Char{\'o}, Gisela D and Ghil, Michael and Sciamarella, Denisse},
  journal={Chaos: An Interdisciplinary Journal of Nonlinear Science},
  volume={33},
  number={10},
  year={2023},
  publisher={AIP Publishing}
}

@InCollection{sciamarella2024new,
  author    = {Sciamarella, Denisse and Char{\'o}, Gisela D.},
  booktitle = {Topological Methods for Delay and Ordinary Differential Equations: With Applications to Continuum Mechanics},
  publisher = {Springer},
  title     = {New elements for a theory of chaos topology},
  year      = {2024},
  pages     = {191--211},
}

@misc{sciamarella2023code,
  author = {Sciamarella, Denisse},
  title = {Templex: A bridge between homologies and templates for chaotic attractors},
  howpublished = {Wolfram Community, STAFF PICKS},
  month = {December},
  year = {2023},
  url = {https://community.wolfram.com/groups/-/m/t/3079776},
  note = {Accessed: 2023-12-07}
}

@article{Svec2025,
  author  = {Šv\v{e}c, Jan G. and Zhang, Zhaoyan},
  title   = {Application of nonlinear dynamics theory to understanding 
             normal and pathologic voices in humans},
  journal = {Philosophical Transactions of the Royal Society B},
  volume  = {380},
  number  = {1923},
  pages   = {20240018},
  year    = {2025},
  doi     = {10.1098/rstb.2024.0018}
}

@Article{Maraldi.ea.2025,
  author    = {Maraldi, B. and Dijkstra, H. A. and Ghil, M.},
  journal   = {Chaos: An Interdisciplinary Journal of Nonlinear Science},
  title     = {Intraseasonal atmospheric variability under climate trends},
  year      = {2025},
  issn      = {1089-7682},
  number    = {5},
  volume    = {35},
  doi       = {10.1063/5.0253103},
  publisher = {AIP Publishing},
}

@article{grandis2003directed,
  title     = {Directed homotopy theory. I},
  author    = {Grandis, Marco},
  journal   = {Cahiers de Topologie et Géométrie Différentielle Catégoriques},
  volume    = {44},
  number    = {4},
  pages     = {281--316},
  year      = {2003}
}

@book{grandis2009directed,
  title     = {Directed Algebraic Topology: Models of Non-Reversible Worlds},
  author    = {Grandis, Marco},
  series    = {New Mathematical Monographs},
  volume    = {13},
  publisher = {Cambridge University Press},
  address   = {Cambridge},
  year      = {2009},
  doi       = {10.1017/CBO9781139193419}
}

@inproceedings{goubault2003some,
  title     = {Some geometric perspectives in concurrency theory},
  author    = {Goubault, Eric},
  booktitle = {Proceedings of the Workshop on Geometry and Topology in Concurrency Theory},
  year      = {2003},
  publisher = {Electronic Notes in Theoretical Computer Science},
  volume    = {81},
  pages     = {1--39},
  doi       = {10.1016/S1571-0661(05)82556-0}
}

@article{gaucher2003homotopical,
  title     = {Homotopical interpretation of globular complex by a model category for homotopy theory of concurrency},
  author    = {Gaucher, Philippe},
  journal   = {Mathematical Structures in Computer Science},
  volume    = {15},
  number    = {3},
  pages     = {409--451},
  year      = {2005},
  doi       = {10.1017/S0960129504004594}
}

@article{birman1983knotted,
  title     = {Knotted periodic orbits in dynamical systems—{I}: {L}orenz’s equations},
  author    = {Birman, Joan S. and Williams, R. F.},
  journal   = {Topology},
  volume    = {22},
  number    = {1},
  pages     = {47--82},
  year      = {1983},
  doi       = {10.1016/0040-9383(83)90049-3}
}

@article{bauer2023cycling,
  title={Cycling signatures: Identifying cycling motions in time series using algebraic topology},
  author={Bauer, Ulrich and Hien, David and Junge, Oliver and Mischaikow, Konstantin},
  journal={arXiv preprint arXiv:2312.04734},
  year={2023}
}

@article{zou2019complex,
  title={Complex network approaches to time series analysis},
  author={Zou, Yong and Donner, Reik V and Marwan, Norbert and Donges, Jonathan F and Kurths, J{\"u}rgen},
  journal={Physics Reports},
  volume={787},
  pages={1--97},
  year={2019},
  publisher={Elsevier}
}

@Book{Philander1990,
  author    = {Philander, S. G. H.},
  publisher = {Academic Press},
  title     = {{ El Ni\~no and the Southern Oscillation.}},
  year      = {1990},
  address   = {New York},
}

\end{document}